\documentclass[11pt]{article}
\usepackage{amsfonts,mathrsfs,bbm,amsmath,latexsym}
\newcommand{\diam}{{\rm diam}}

\newcommand{\supp}{{\rm supp}}

\newcommand{\rmd}{\mbox{\rm d}}
\newcommand{\rep}{\mbox{\rm Re}}

\newcommand{\rme}{\mbox{\rm e}}

\newcommand{\bfA}{{\mbox{\boldmath $A$}}}

\title{Mathematical analysis of population migration and its effects to spread of
  epidemics\footnote{This work is supported by the China National Natural Science Fundation
  under grant number 11171357.}}
\author{Shangbin Cui$^1$\footnote{Corresponding author. E-mail:  cuishb@mail.sysu.edu.cn}\;\;\;\;
  and\;\; Meng Bai$^2$}
\date{{\small 1: Department of Mathematics, Sun Yat-Sen University,
   Guangzhou, Guangdong 510275, China} \\ [-0.05cm]
  {\small 2: School of Mathematics and Statistics, Zhaoqing University,
   Zhaoqing, Guangdong 526061, China}}

\begin{document}
\maketitle
\begin{abstract}
  In this paper we study some mathematical models describing evolution of population density and
  spread of epidemics in population systems in which spatial movement of individuals depends only
  on the departure and arrival locations and does not have apparent connection with the population
  density. We call such models as population migration models and migration epidemics models,
  respectively. We first apply the theories of positive operators and positive semigroups to make
  systematic investigation to asymptotic behavior of solutions of the population migration models
  as time goes to infinity, and next use such results to study asymptotic behavior of solutions of
  the migration epidemics models as time goes to infinity. Some interesting properties of solutions
  of these models are obtained.



\vspace*{0.5cm}

\centerline{\bf Table of Contents}
\begin{center}
\begin{minipage}[t]{12cm}
\begin{flushleft}
\begin{description}
\item[]\S1 Introduction
\item[]\S2 Population migration models
\begin{description}
\item[]\hspace*{-0.7cm} 2.1 The migration rate function
\item[]\hspace*{-0.7cm} 2.2 The proliferation-stationary population migration model
\item[]\hspace*{-0.7cm} 2.3 Solutions of the equation (2.3) in some simple situations
\item[]\hspace*{-0.7cm} 2.4 Asymptotic behavior of solutions in the ergodic migration case
\item[]\hspace*{-0.7cm} 2.5 The proliferation non-stationary population migration model
\item[]\hspace*{-0.7cm} 2.6 A short discussion to non-ergodic cases
\end{description}
\item[]\S3 Migration epidemics models
\begin{description}
\item[]\hspace*{-0.7cm} 3.1 The migration SI model
\item[]\hspace*{-0.7cm} 3.2 The migration SIR model
\item[]\hspace*{-0.7cm} 3.3 The migration SIRE model
\end{description}
\item[]\S4 Conclusions
\end{description}
\end{flushleft}
\end{minipage}
\end{center}
\end{abstract}
\section{Introduction}
\setcounter{equation}{0}

\hskip 2em
  Since publications of the pioneering work of Ross \cite{Ross}, McKendrick \cite{McK} and Kermack
  and McKendrick \cite{KerMcK} on mathematical analysis of dynamics of infectious diseases, a vast
  variety of mathematical models aiming at revealing the mechanisms of spread of epidemics have
  appeared in the literature, and great progress has been made in investigation of the dynamics of
  infectious diseases by using such models \cite{Cap}. Classical epidemics models only consider
  temporal effects of interactions of different groups of population such as the susceptible group,
  the infective group, and the recovered (or removed or immune) group, and spatial distribution of
  the population in which epidemics is spread are usually ignored. Since spatial movement of population
  is clearly an important factor of the spread of epidemics, recently increasingly much attention
  on this topic has been paid on investigation of spatial effects of the population on dynamical
  behavior of infectious diseases, and an increasingly large amount of spatial epidemics models
  have appeared in the literature, cf. \cite{ABLN1,ABLN2,ADH,AD1,AD2,Cos3,FLM,GLW,GW,GWB,JW,LM,
  Mur,PL,Ruan,SD,WM,WZ} and references therein.

  A basic concern of all spatial epidemics models is how the host population of the epidemics is
  distributed, moved and developed in the habitat domain. It follows that applications of
  different temporal-spatial population models to the study of dynamics of infectious diseases
  yield different spatial epidemics models. Existing spatial epidemics models are either built
  on the basis of multi-patch (or multi-country, metapopulation etc.) population models
  \cite{ABLN1,ADH,AD1,AD2,Cos3,JW,LM,SD,WM,WZ}, or constructed on applications of reaction-diffusion
  equations to population systems \cite{ABLN2,FLM,GLW,GW,GWB,Mur,PL,Ruan}, or obtained by combining
  ideas of the above two approaches together. Mathematical analysis of such models have revealed
  many interesting phenomena, cf. the references mentioned above, especially
  \cite{ABLN1,ABLN2,AD2,Cos3,GLW,Ruan}.

  In the meantime of sufficiently recognizing great success achieved by those spatial epidemics
  models as mentioned above, we must also be aware of their limitations in application to the
  spread of epidemics in the human population system. Clearly, multi-patch epidemics models are
  useful only to population systems with habitat domain being a combination of isolated
  subdomains like pools, lakes and islands. Moreover, lakes and islands must be very small so that
  when an infectious disease occurs in one lake or island then all individuals living in that
  lake or island are infected within a short period, i.e., population in each patch must be
  homogeneous in getting infected by the epidemics. It follows that application of such models to
  continental human population system is not very suitable. As far as reaction-diffusion epidemics
  models are concerned, its limitation in application arises from a different argument. In ancient
  times transportation tools are so undeveloped that common people usually cannot move very far
  from their habitat locations. Taking this situation into account, reaction-diffusion equations
  are a good choice to mimic spatial movement of ancient human population in which epidemics is
  spread. This explains success of application of the reaction-diffusion epidemic model in the
  analysis of the spread of black death in Europe during 1347$\sim$1350 \cite{Mur}.  But nowadays
  due to popular applications of long-distance reachable high-speed transportation tools like
  trains and airplanes, travel distance of common people has dramatically increased. It turns out
  that the situation has become quite different, and many non-diffusion types of spatial movement
  of population such as tourist and commercial travels, migrations of rural people to cities for
  study and jobs, visiting travels between relative families or friends and so on have become main
  causes of spread of infectious disease. Take the worldwide spread of severe acute respiratory
  syndrome (SARS) during November 2002 till July 2003 as an example. It first appeared in November
  16, 2002 in Foshan, which is a county city in the Guangdong Province adjacent to Guangzhou.
  Weeks after its onset, some patients were sent to Guangzhou for medical treatment, resulting in
  the spread of the disease in hospitals in Guangzhou during the following three months. Infected
  people during that period were mainly medical personnel and patients in hospitals in Guangzhou,
  and their families. On February 21, 2003, a retired professor of the Second Affiliated Hospital
  of the Sun Yat-Sen University who had been keeping working in his hospital as a doctor went to
  Hongkong to attend the wedding ceremony of a relative of him, bringing SARS to Hongkong. It is
  in Hongkong that spread of SARS first became explosive: 1755 people were infected and among them
  299 died during the following four months. On March 6, 2003 a hospital in Beijing received the
  first SARS patient. Soon Beijing became the second explosive spread place of this disease: 2434
  people were infected and among them 147 died during the following four months. On March 13, 2003
  Taiwan reported first SARS infection case, and during the following four months Taiwan totally
  reported 307 infected cases (including 47 deaths). According to a report made by WHO, from
  November 16, 2002 till July 11, 2003 totally 29 countries in the world reported SARS infection
  cases, and totaly 8069 people were infected by SARS with 775 people among them died \cite{SARS1,
  SARS2}. Besides the cities and regions mentioned above, the other very serious SARS infection
  places include the following provinces of China: Guangdong (1514 cases with 56 deaths), Shanxi
  (445 cases with 20 deaths), Inner Mongolia (289 cases with 25 deaths), Hebei (210 cases with 10
  deaths) and Tianjin (176 cases with 12 deaths) (see Table 1 for the distribution in China of
  SARS cases reported in the 2002$\sim$2003 outbreak), and the following countries: Canada (250
  cases with 38 deaths), Singapore (206 cases with 32 deaths), United States (71 cases without
  deaths) and Viet Nam (63 cases with 5 deaths) (see Table 2 for the distribution of SARS cases
  reported in the 2002$\sim$2003 outbreak). Here we particularly mention that though in the
  Guangdong Province (whose capital city is Guangzhou) the spread of SARS is serious, the situation
  in its four neighboring provinces is, however, much better: There were only 22 cases (with 3
  deaths) in the Guangxi Province, 6 cases (with 1 death) in the Hunan Province, 1 case (without
  death) in the Jiangxi Province, and 3 cases (without death) in the Fujian Province. From this
  example we easily see that in modern world the spread of epidemics among human population cannot
  be well modeled by reaction-diffusion equations, but instead seriously depends on communications
  of inhabitants in different places.

  Based on the above consideration, in this paper we study a different class of spatial epidemics
  models constructed on the basis of population migration models. Here, population migration models
  refer to so-called {\em population dispersal models} in existing literatures \cite{Bat,Can2,Can3,
  Cos4,Hut2,Cov1,Cov2,Cov3,Cov4}. We use the terminology {\em migration} rather than {\em dispersal}
  due to the following consideration: These models are constructed on the assumption that spatial
  movement of individuals in the population system under study depends only on geographic variables,
  i.e., the tendency and strength of population movement from one location the the other is uniquely
  determined by the departure and the arrival locations (because of some population-density
  independent reasons such as pursuit of better resources in job and education, official, commercial
  or sightseeing-aimed travel, and etc.); it does not have apparent dependence on the density of the
  population, so that both ``diffusion'' and ``dispersal'' are not suitable to be used to depict
  such a movement. Hence, these population models are not new; they have appeared many times in the
  above-mentioned and some other references. However, to the best of our knowledge, dynamical
  behavior of the population system described by these models has not been very well understood.
  For the purpose to get a good understanding of dynamics of the spread of infectious disease in
  such a population system, we shall thus first make a systematic investigation to asymptotic
  behavior of solutions of the population migration models. This will be performed by using the
  well-established theories of positive operators and positive semigroups as lectured in the
  references \cite{Nag,EN,Neer,Sch}, and will be presented in Section 2. We shall next in Section 3
  use the results obtained from a such analysis to study three migration epidemics models: The
  migration SI model, the migration SIR model, and the migration SIRE model. These temporal-spatial
  epidemics models are obtained by adding migration terms into classical spatial homogeneous
  epidemics models; they have the form of nonlinear differential-integral equations. As a first
  step of the study in this direction, we shall be content with deriving some elementary results on
  asymptotic behavior of solutions of these migration epidemic models, and leave more advanced
  investigation of other related topics, which are very diversified, for future work. In the last
  section we make a summary to the discussion of this manuscript.

  {\bf Notations}\ \ For reader's convenience we collect some notations to be used later as follows:
\begin{itemize}
\item $\Omega$ represents a bounded domain in $\mathbb{R}^n$. $|\Omega|$ denotes Lebesque measure
  of $\Omega$.
\item $C(\overline{\Omega})$ and $L^1(\Omega)$ denote the Banach spaces of all continuous functions
  in $\overline{\Omega}$ and all Lebesque integrable functions in $\Omega$, respectively, with norms
\begin{equation*}
  \|\varphi\|_{\infty}=\displaystyle\max_{x\in\overline{\Omega}}|\varphi(x)|\;\;
  (\mbox{for}\;\;\varphi\in C(\overline{\Omega})) \quad \mbox{and} \quad
  \|\varphi\|_1=\displaystyle\int_{\Omega}|\varphi(x)|\rmd x\;\;
  (\mbox{for}\;\;\varphi\in L^1(\Omega))
\end{equation*}
  respectively.
\item For $\varphi\in C(\overline{\Omega})$ (resp. $\varphi\in L^1(\Omega)$), $\varphi\geq 0$ means
  that $\varphi(x)\geq 0$ for all $x\in\Omega$ (resp. for almost every $x\in\Omega$), $\varphi>0$
  means that $\varphi\geq 0$ and $\varphi\neq 0$, and $\varphi\gg 0$ means that $\varphi(x)>0$ for
  all $x\in\overline{\Omega}$ (resp. for almost every $x\in\Omega$). If $\varphi\geq 0$ then
  $\varphi$ is said to be {\em positive} in $\Omega$, and if $\varphi\gg 0$ then $\varphi$ is said
  to be {\em strictly positive} in $\Omega$. For $\varphi,\psi\in C(\overline{\Omega})$ or
  $L^1(\Omega)$), $\varphi\geq\psi$ means that $\varphi-\psi\geq 0$, and $\varphi\leq\psi$ means
  that $\psi\geq\varphi$. It is well-known that ordered in this way, $C(\overline{\Omega})$ and
  $L^1(\Omega)$ are Banach Lattices (cf. \cite{Nag,Sch}).
\item $\mathcal{M}(\Omega)$ denotes the Banach space of all bounded regular Borel measures on
  $\Omega$, with norm $\|\rmd\mu\|=\displaystyle\int_{\Omega}|\rmd\mu|=\int_{\Omega}\rmd\mu^+
  +\int_{\Omega}\rmd\mu^-$ for $\rmd\mu\in\mathcal{M}(\Omega)$, where $\rmd\mu^+$, $\rmd\mu^-$
  denote the positive and negative parts of $\rmd\mu$, respectively. The well-known Riesz
  representation theorem asserts that $(C(\overline{\Omega}))'=\mathcal{M}(\Omega)$, i.e., for
  any $f\in (C(\overline{\Omega}))'$ there exists an unique $\rmd\mu\in\mathcal{M}(\Omega)$ such
  that
\begin{equation*}
  \langle\varphi,f\rangle=f(\varphi)=\displaystyle\int_{\Omega}\varphi(x)\rmd\mu(x)=
  \int_{\Omega}\varphi(x)\rmd\mu^+(x)-\int_{\Omega}\varphi(x)\rmd\mu^-(x)
\end{equation*}
  for all $\varphi\in C(\overline{\Omega})$, and $\|f\|=\|\rmd\mu\|$. The notation $\rmd\mu\geq 0$
  means that $\rmd\mu$ is a positive measure. It follows that $\mathcal{M}(\Omega)$ is also a
  Banach lattice. Moreover, $\rmd\mu>0$ means that $\rmd\mu\geq 0$ and $\rmd\mu\neq 0$, and $\rmd
  \mu\gg 0$ means that $\rmd\mu(O)>0$ for any nonempty open subset $O$ of $\Omega$, or equivalently,
  $\displaystyle\int_{\Omega}\varphi(x)\rmd\mu(x)>0$ for any $\varphi\in C(\overline{\Omega})$ such
  that $\varphi>0$ (cf. \cite{Nag,Sch}).
\item For $\varphi\in C(\overline{\Omega})$, $\supp\varphi$ denotes the support of $\varphi$,
  i.e., $\supp\varphi=$ closure of the set $\{x\in\overline{\Omega}:\varphi(x)>0\}$. For
  $\varphi\in L^1(\Omega)$, $\supp\varphi$ denotes the essential support of $\varphi$, i.e.,
  $\supp\varphi=$ closure (in $\Omega$) of the set $\{x_0\in\Omega:\mbox{there is a neighborhood
  $O$ of $x_0$ such that}\;\;\varphi(x)>0\;\;\mbox{for a. e. }\;\,x\in O\}$. For $\rmd\mu\in
  \mathcal{M}(\Omega)$, $\supp\rmd\mu$ denotes the support of $\rmd\mu$, i.e., $\supp\rmd\mu=$
  closure of the largest open subset $O$ of $\Omega$ such that $\rmd\mu|_O\gg 0$.
\item For a Banach space $X$, $\mathscr{L}(X)$ denotes the Banach algebra of all bounded linear
  operators on $X$, and $I$ denotes the unit operator or the identity map in $X$. In case we need
  to emphasize that it is the identity map in the space $X$, we use $id_X$ to replace $I$. If
  $A\in\mathscr{L}(X)$ then its operator norm is denoted as $\|A\|_{\mathscr{L}(X)}$.
\item For $A\in\mathscr{L}(X)$, $\rho(A)$ and $\sigma(A)$ denote the resolvent set and spectrum
  of $A$, respectively. For $\lambda\in\rho(A)$, $R(\lambda,A)$ denotes the resolvent of $A$ (at
  the point $\lambda$), i.e., $R(\lambda,A)=(\lambda I-A)^{-1}$.
\item For $A\in\mathscr{L}(X)$, $r(A)$ and $s(A)$ denote the spectral radius and the spectral bound
of $A$, respectively, i.e.,
\begin{equation*}
  r(A)=\sup\{|\lambda|:\lambda\in\sigma(A)\}, \quad
  s(A)=\sup\{\rep\lambda:\lambda\in\sigma(A)\}.
\end{equation*}
\item For $A\in\mathscr{L}(X)$, $\rme^{tA}$ ($t\geq 0$) denotes the uniformly continuous semigroup
  in $X$ generated by $A$, i.e., $\rme^{tA}=\displaystyle\sum_{k=0}^{\infty}\frac{t^k}{k!}A^k$, and
  $\omega_0(A)$ denotes the growth bound of the semigroup $\rme^{tA}$, i.e.,
\begin{equation*}
  \omega_0(A)=\inf\{\omega\in\mathbb{R}:\exists M>0\;\,\mbox{s.t.}\;\,\|\rme^{tA}\|_{\mathscr{L}(X)}
  \leq Me^{\omega t},\forall t\geq 0\}.
\end{equation*}
  Recall that the spectral mapping theorem (cf. Lemma 3.13 in Chapter I of \cite{EN}) says that
\begin{equation*}
  \sigma\big(\rme^{tA}\big)=\rme^{t\sigma(A)}:=\{\rme^{t\lambda}:\lambda\in\sigma(A)\},
  \;\;\forall t\geq 0,
\end{equation*}
  so that $\omega_0(A)=s(A)$.
\item Finally, for $A\in\mathscr{L}(X)$, $\sigma_{ess}(A)$ and $\omega_{ess}(A)$ denote the essential
  spectrum of $A$ and the essential growth bound of the semigroup $\rme^{tA}$, respectively, cf.
  Discussion 1.20 and Definition 2.9 in Chapter IV of \cite{EN}, respectively, and also $\S3.6$ of
  \cite{Neer},  for definitions.
\end{itemize}

\section{Population migration models}
\setcounter{equation}{0}

\hskip 2em
  In this section we make a systematic study to some population migration models. These models have
  the form of differential-integral equations. They describe the evolution of the density of a
  system of population living in a spatial domain $\Omega\subseteq\mathbb{R}^n$. Individuals in
  this population system can move from one location to the other with a rate uniquely determined by
  their departure and arrival locations. We use the terminology ``migration'' to address spatial
  movement of the population and call the rate mentioned above as ``migration rate''. The migration
  rate function determines to a great extent the dynamical behavior of the population system. We
  shall mainly consider a special class of population system in which individuals can migrate from
  any one location to any other location through finite steps of intermediate migrations. We call
  such a population system as ``ergodic system''. In the proliferation-stationary case, dynamical
  behavior of the population in such a system can be very well understood by applying the theories
  of irreducible positive operators and irreducible positive semigroups. We next simply consider
  two other classes of population systems, one of which has the property that the whole system can
  be divided into several mutually independent subsystems, and the other of which has the property
  that the whole system can be divided into two subsystems, with population in one of these two
  subsystems being able to migrate to the other whereas oppositely unable. From analysis to these
  simple examples we can get a rough impression on the dynamical behavior of the other more general
  population systems. In the last part of this section we study dynamics of ergodic population
  system in the case that its proliferation is not stationary.

  The results obtained in this section will be used to study dynamical behavior of solutions of
  some migration epidemics models in the next section.

\subsection{The migration rate function}

\hskip 2em
  Let $\Omega$ be a domain in $\mathbb{R}^N$ in which the population lives. We introduce
  a {\em migration rate function} $v(x,y)$ in $\Omega$, or more precisely, the {\em migration rate}
  from location $y\in\Omega$ to location $x\in\Omega$, defined as follows: For an off-diagonal point
  $(x,y)\in\Omega\times\Omega$, $x\neq y$, let $O_x$ and $O_y$ be two disjoint sub-domains of $\Omega$
  containing $x$ and $y$, respectively. Assume that within an unit time interval the total amount
  of population migrated from $O_y$ to $O_x$ is $M(O_x,O_y)$ and the total amount of population in
  $O_y$ before migration is $M(O_y)$. Then
\begin{equation*}
  v(x,y)=\lim_{\diam(O_x)\to 0\atop\diam(O_y)\to 0}\frac{M(O_x,O_y)}{|O_x|M(O_y)},
\end{equation*}
  where $|O_x|$ denotes the measure of the domain $O_x$. For a diagonal point of $\Omega\times\Omega$
  we define the value of $v(x,y)$ to be its off-diagonal limit. Let
\begin{equation}
  v_i(x)=\int_{\Omega}v(x,y){\rm d}y, \qquad   v_e(y)=\int_{\Omega}v(x,y){\rm d}x;
\end{equation}
  they are called the {\em immigration rate} and {\em emigration rate} at the locations
  $x$ and $y$, respectively. Notice that these definitions immediately yield the following
  relation:
\begin{equation}
  \int_{\Omega}v_i(x){\rm d}x=\int_{\Omega}v_e(y){\rm d}y.
\end{equation}

  From the definition of the migration rate function we see that for a point $x\in\Omega$, if
  $v(x,y)=0$ for all $y\in\Omega$, then there is no migration from other location to the
  location $x$; in this case $v_i(x)=0$. For a point $x\in\Omega$, if $v(y,x)=0$ for all
  $y\in\Omega$, then there is no migration from the location $x$ to other location; in
  this case $v_e(x)=0$. If at a location $x$ there holds $v_i(x)>v_e(x)$ then in a small
  neighborhood of $x$ the amount of populations immigrated in is larger than the amount of
  populations emigrated out, whereas if $v_e(x)>v_i(x)$ then in a small neighborhood of $x$
  the amount of populations emigrated out is larger than the amount of populations immigrated
  in. In the case $v_e(x)=v_i(x)$ we say that immigration and emigration at the location $x$
  is balanced. In particular, if the relation $v_e(x)=v_i(x)$ holds for all
  $x\in\Omega$, we call such a migration as {\em balanced migration}. A particular case
  of the balanced migration is {\em symmetric migration}, which means that $v(x,y)=v(y,x)$
  for all $x,y\in\Omega$.

  From the relation (2.2) we see that if there is a subdomain $\Omega_1\subseteq\Omega$
  in which $v_e(x)>v_i(x)$ then there must be another subdomain $\Omega_2\subseteq\Omega$
  in which $v_e(x)<v_i(x)$, and vice versa.

  Dynamical behavior of a population system depends in great extent on the location of
  $\supp v$ in $\Omega\times\Omega$. For instance, if $\Omega$ is divided into two disjoint
  parts $\Omega_1$ and $\Omega_2$, i.e., $\Omega_1\cap\Omega_2=\emptyset$ and $\Omega=\Omega_1
  \cup\Omega_2$, such that $v(x,y)=0$ whenever $x$ and $y$ are not in the same part, i.e.,
  $v(x,y)=0$ if either $(x,y)\in\Omega_1\times\Omega_2$ or $(x,y)\in\Omega_2\times\Omega_1$,
  then clearly the whole population system can be divided into two independent subsystems defined
  in $\Omega_1$ and $\Omega_2$ respectively, with migration rate functions $v|_{\Omega_1\times
  \Omega_1}$ and $v|_{\Omega_2\times\Omega_2}$ respectively; the development of population in
  each subsystem does not have any influence to that of the other subsystem. On the other hand,
  as we shall see later, if $v(x,y)>0$ for all $(x,y)\in\Omega\times\Omega$ then the development
  of population in the whole domain $\Omega$ is integrative: if the initial population density
  $N_0(x)$ is not identically vanishing in $\Omega$ then the population will eventually spread
  to the whole domain $\Omega$ with an everywhere nonvanishing distribution $N_{\infty}(x)$
  determined by the migration rate function $v$ and the total amount of the initial population
  $M_0=\displaystyle\int_{\Omega}N_0(x)\rmd x$, regardless how the initial population is
  distributed in the domain $\Omega$ (i.e., $N_{\infty}(x)$ does not depend on specific profile
  of the initial function $N_0(x)$; in particular, if $N_0(x)>0$ in a small neighborhood of a
  point $x_0\in\Omega$, then $N_{\infty}(x)>0$ for all $x\in\Omega$, even if $N_0(x)$ vanishes
  outside that neighborhood of $x_0$). Another interesting example is that as before $\Omega=
  \Omega_1\cup\Omega_2$ with $\Omega_1\cap\Omega_2=\emptyset$, but now $v(x,y)>0$ for $(x,y)\in
  \Omega_1\times\Omega_1$, $\Omega_1\times\Omega_2$ and $\Omega_2\times\Omega_2$, and $v(x,y)=0$
  for $(x,y)\in\Omega_2\times\Omega_1$. In this example, the whole system can be divided into
  two mutually dependent subsystems, one defined in the subdomain $\Omega_1$ and the other defined
  in the subdomain $\Omega_2$. Clearly, dependence of the two subsystems is biased: Since there
  exists positive migration from the subsystem defined in the subdomain $\Omega_2$ to that in the
  subdomain $\Omega_1$ but there does not exist opposite-direction migration, development of
  population in the subdomain $\Omega_2$ has great influence to that in the subdomain $\Omega_1$,
  but development of population in the subdomain $\Omega_1$ does not have any effect to that in
  the subdomain $\Omega_2$. See the subsections 2.4--2.6 for proofs of these assertions.

  The above simple examples also show that dynamical behavior of a population system with
  migration is very complex. In this paper we do not pursue a complete understanding to an
  arbitrary population system with migration. Instead, we mainly consider dynamical behavior of
  a special class of such system as defined below:
\medskip

  {\bf Definition 2.1}\ \ {\em We say that a population system is {\em ergodic} if its migration
  rate function $v$ satisfies the following property: For any point $(x_0,y_0)\in\Omega\times
  \Omega$, $x_0\neq y_0$, there exist finite number of points $(x_1,y_1),(x_2,y_2),\cdots,
  (x_m,y_m)\in\Omega\times\Omega$ such that
$$
  y_1=y_0, \;\; y_2=x_1,  \;\; y_3=x_2, \;\; \cdots,  \;\; y_m=x_{m-1}, \;\; x_m=x_0,
$$
  and for each $1\leq j\leq m$, $v(x,y)>0$ for all $(x,y)$ in a neighborhood of $(x_j,y_j)$.
  In this case we also say that the function $v$ is ergodic.}
\medskip

  From the above definition we see that if a population system is ergodic then it is possible
  for an individual in this system at any location to migrate to any other location through
  finite steps of intermediate migrations. Clearly, for such a population system we have
$$
  v_i(x)>0 \;\; \mbox{and} \;\; v_e(x)>0 \quad \mbox{for all}\;\; x\in\Omega.
$$
  It is also clear that if $v(x,y)>0$ for all $(x,y)\in\Omega\times\Omega$ then the population
  system is an ergodic system; in this case we say that the population system is {\em completely
  ergodic}. An example of a population system which is ergodic but not completely ergodic is
  that the migration rate function $v$ satisfies the following {\em Coville condition} \cite{Cov1}
  $\sim$ \cite{Cov4}:
$$
  \mbox{\em There exist}\;\;\varepsilon_0>0\;\;\mbox{\em and}\;\;c_0>0\;\;\mbox{\em such that}\;\;
  \inf_{x\in\Omega}\Big(\inf_{y\in B(x,\varepsilon_0)}v(x,y)\Big)\geq c_0.
$$
  Another example is as follows: Let $\Omega=[0,1]$ and $v$ be a nonnegative function defined in
  $[0,1]\times [0,1]$ having the following property:
$$
  \supp v=\Big\{(x,y)\in [0,1]\times [0,1]:\frac{1}{4}\leq x+y\leq\frac{3}{4}\;\;
  \mbox{or}\;\;\frac{5}{4}\leq x+y\leq\frac{7}{4}\Big\}.
$$
  It is easy to check that the population system with $v$ being its migration rate function is
  ergodic but not completely ergodic and it also does not satisfies the Coville condition.

  Recall that a positive bounded linear operator $T$ in a Banach lattice $X$ is said to be {\em
  irreducible} if for any $0<x\in X$ (i.e. $x\geq 0$ and $x\neq 0$) and $0<f\in X'$ (i.e. $f\geq 0$
  and $f\neq 0$) there exists corresponding positive integer $m$ such that $\langle T^mx,f\rangle
  >0$, cf. Section 8 in Chapter III of \cite{Sch}. Recall also that a positive $C_0$-semigroup
  $T(t)$ ($t\geq 0$) in a Banach lattice $X$ is said to be  {\em irreducible} if for any $0<x\in X$
  and $0<f\in X'$ there exists corresponding $t_0>0$ such that $\langle T(t_0)x,f\rangle>0$, cf.
  Section 3 in Chapter C-III of \cite{Nag}. Let $v$ be a nonnegative continuous function defined in
  $\overline{\Omega}\times\overline{\Omega}$, where $\Omega$ is a bounded domain in $\mathbb{R}^n$.
  Let $X$ be any one of the function spaces $C(\overline{\Omega})$ and $L^p(\Omega)$, $1\leq p<
  \infty$. Let $K$ be the integral operator in $X$ with kernel $v$, i.e., $K\varphi(x)=
  \displaystyle\int_{\Omega}v(x,y)\varphi(y)\rmd y$ for $\varphi\in X$. For any real-valued
  continuous function $b$ in $\overline{\Omega}$, let $M_b$ denote the multiplication operator in
  $X$ related to the function $b$, i.e., $M_b\varphi(x)=b(x)\varphi(x)$ for $\varphi\in X$.
\medskip

  {\bf Lemma 2.2}\ \ {\em Let notations be as above and assume that $v$ is ergodic. We have the
  following assertions:
  $(1)$ For any nonnegative continuous function $b$ in $\overline{\Omega}$, the positive operator
  $K+M_b$ in $X$ is irreducible.
  $(2)$ For any continuous function $b$ in $\overline{\Omega}$, the positive uniformly continuous
  semigroup $e^{t(K+M_b)}$ in $X$ generated by $K+M_b$ is irreducible.}
\medskip

  {\em Proof}:\ \ We give a proof only for the case $X=C(\overline{\Omega})$, because argument for
  the other case $X=L^p(\Omega)$ ($1\leq p<\infty$) is similar. Let $0<\varphi\in
  C(\overline{\Omega})$ and $0<f=\rmd\mu\in\mathcal{M}(\Omega)=(C(\overline{\Omega}))'$. Let $x_0,
  y_0\in\Omega$ be such that $\rmd\mu(y_0)>0$ (i.e., $\displaystyle\int_{\Omega}\varphi(x)\rmd
  \mu(x)>0$ for any $0\leq\varphi\in C(\overline{\Omega})$ such that $\varphi(y_0)>0$) and
  $\varphi>>0$ in some neighborhood $\Omega_1$ of $y_0$. By ergodicity of the function $v$, there
  exist finite number of points $(x_1,y_1),(x_2,y_2),\cdots,(x_m,y_m)\in\Omega\times\Omega$ such
  that
$$
  y_1=y_0, \;\; y_2=x_1,  \;\; y_3=x_2, \;\; \cdots,  \;\; y_m=x_{m-1}, \;\; x_m=x_0,
$$
  and for each $1\leq j\leq m$, $v(x,y)>0$ for all $(x,y)$ in a neighborhood of $(x_j,y_j)$. It
  follows that for each $2\leq j\leq m$ there exists a corresponding neighborhood $\Omega_j$ of
  $y_j=x_{j-1}$ such that $v(x,y)>0$ for all $(x,y)\in\Omega_j\times\Omega_{j-1}$ ($\Omega_1$
  probably has to be replaced with a smaller one when necessary). Hence
\begin{eqnarray*}
  \langle K^m\varphi,f\rangle
  &=&\displaystyle\int_{\Omega}\int_{\Omega}\cdots\int_{\Omega}\int_{\Omega}
  v(x,\xi_m)\cdots v(\xi_3,\xi_2)v(\xi_2,\xi_1)\varphi(\xi_1)
  \rmd \xi_1\rmd \xi_2\cdots\rmd \xi_m\rmd\mu(x)
\\
  &\geq &\displaystyle\int_{\Omega}\int_{\Omega_m}\cdots\int_{\Omega_2}\int_{\Omega_1}
  v(x,\xi_m)\cdots v(\xi_3,\xi_2)v(\xi_2,\xi_1)\varphi(\xi_1)
  \rmd \xi_1\rmd \xi_2\cdots\rmd \xi_m\rmd\mu(x)
\\
  &>& 0.
\end{eqnarray*}
  This shows that the integral operator $K$ is irreducible. From this assertion we easily see that
  for any nonnegative continuous function $b$ in $\overline{\Omega}$, the positive operator
  $K+M_b$ is also irreducible. This proves the assertion (1). Next, since $\rme^{tK}=\displaystyle
  \sum_{m=0}^{\infty}\frac{t^m}{m!}K^m$, from the above assertion we also easily see that the
  semigroup $\rme^{tK}$ is irreducible. By Proposition 3.3 in Chapter B-III of \cite{Nag}, it
  follows that for any continuous function $b$ in $\overline{\Omega}$, the semigroup $e^{t(K+M_b)}$
  is also irreducible. This proves the lemma. $\quad\Box$

\subsection{Proliferation-stationary population migration model}

\hskip 2em
  We denote by $N(x,t)$ the {\em population density} at the location $x\in\Omega$ at time $t$. We
  assume that the population at every location in $\Omega$ is in the proliferation-stationary
  state, i.e., the birth rate and the death rate are equal at every point $x\in\Omega$, and only
  consider the effect of migration of population. For an arbitrary sub-domain $O\subseteq\Omega$
  and an arbitrary time interval $[t_1,t_2]$, we have
\begin{eqnarray*}
  \underbrace{\int_{O}N(x,t_2)\rmd x-\int_{O}N(x,t_1)\rmd x}_{\mbox{\rm increment of population}}
  &\;\;=\;\;&\underbrace{\int_{t_1}^{t_2}\int_{O}\int_{\Omega\backslash O}v(x,y)N(y,t)
  \rmd y\rmd x\rmd t}_{\mbox{\rm population immigrated in}}
\\
  &&-\underbrace{\int_{t_1}^{t_2}\int_{\Omega\backslash O}\int_{O}v(x,y)N(y,t)
  \rmd y\rmd x\rmd t}_{\mbox{\rm population emigrated out}}.
\end{eqnarray*}
  Dividing both sides with $(t_2-t_1)|O|$ and next letting $\diam(O)\to 0$ and $t_2-t_1\to 0$,
  we obtain the following differential-integral equation:
\begin{equation}
  \frac{\partial N(x,t)}{\partial t}=\int_{\Omega}v(x,y)N(y,t)\rmd y-v_e(x)N(x,t), \quad
  x\in\Omega,\;\; t>0.
\end{equation}
  We call it the {\em proliferation-stationary population migration equation}, or simply
  {\em population migration equation} later on. We impose the following initial value condition:
\begin{equation}
  N(x,0)=N_0(x), \quad x\in\Omega,
\end{equation}
  where $N_0$ is a given nonnegative function.

  For simplicity we only consider the case where $\Omega$ is a bounded domain, though much of the
  discussion also holds for unbounded domains. Hence later on we always assume that $\Omega$ is a
  bounded domain. We have the following easy and basic assertions:

  (1)\ \ Since the equation (2.4) is a linear differential-integral equation, by using either the
  standard Picard iteration method or the uniformly continuous semigroup theory, we can easily
  prove that under suitable assumptions on the migration function $v$, the initial value
  problem (2.3)--(2.4) is globally well-posed in the function spaces $C(\overline{\Omega})$ and
  $L^p(\Omega)$ for $1\leq p<\infty$. More precisely, if $v\in C(\overline{\Omega}\times
  \overline{\Omega})$ then for any $N_0\in C(\overline{\Omega})$ the problem (2.3)--(2.4) has an
  unique solution $N\in C(\overline{\Omega}\times[0,\infty))$, and the map $N_0\mapsto N$ from
  $C(\overline{\Omega})$ to $C(\overline{\Omega}\times[0,T])$ is linear and continuous for any
  $T>0$. Similarly, if $v\in L^{\infty}(\Omega\times\Omega)$ then for any $1\leq p<\infty$ and $N_0
  \in L^p(\Omega)$ the problem (2.3)--(2.4) has an unique solution $N\in C([0,\infty),L^p(\Omega))$,
  and the map $N_0\mapsto N$ from $L^p(\Omega)$ to $C([0,T],L^p(\Omega))$ is linear and continuous
  for any $T>0$.

  (2)\ \ If $N_0(x)\geq 0$ for all $x\in\Omega$, then  $N(x,t)\geq 0$ for all $x\in\Omega$
  and $t\geq 0$.

  The proofs of the above two assertions are standard, so that we omit them here.

  (3)\ \ The total amount of the population is constant, i.e., letting $M(t)=\displaystyle
  \int_{\Omega}N(x,t)\rmd x$ be the total amount of the population at time $t$ and $M_0=
  \displaystyle\int_{\Omega}N_0(x)\rmd x$ be the initial total amount of the population, we have
\begin{equation}
  M(t)=M_0 \quad \mbox{for all}\;\; t\geq 0.
\end{equation}
  Indeed, since $\displaystyle\int_{\Omega}v(x,y)\rmd x=v_e(y)$, by integrating both sides of
  the equation (2.4) with respect to the variable $x$, we get
$$
  \frac{\rmd M(t)}{\rmd t}=\int_{\Omega}v_e(y)N(y,t)\rmd y-\int_{\Omega}v_e(x)N(x,t)\rmd x=0
   \quad \mbox{for all}\;\; t>0,
$$
  so that $M(t)=M_0$ for all $t\geq 0$.

  (4)\ \ If $v_e(x_0)=v_i(x_0)=0$ for some $x_0\in\Omega$, then $N(x_0,t)=N_0(x_0)$ for all
  $t\geq 0$. Indeed, since $v(x,y)\geq 0$ for all $x,y\in\Omega$, the condition $v_i(x_0)=
  \displaystyle\int_{\Omega}v(x_0,y){\rm d}y=0$ implies that $v(x_0,y)=0$ for all $y\in\Omega$.
  Hence, from the condition $v_e(x_0)=v_i(x_0)=0$ we see that at the point $x_0$ the equation
  (2.4) takes the form
$$
  \frac{\partial N(x_0,t)}{\partial t}=0 \quad \mbox{for}\;\; t>0,
$$
  so that $N(x_0,t)=N_0(x_0)$ for all $t\geq 0$. This means that if at a location the population
  does neither migrate in nor migrate out, then the population density keeps constant at that
  location.

  (5)\ \ If $v_i(x_0)=0$ and $v_e(x_0)>0$ for some $x_0\in\Omega$, then $\displaystyle
  \lim_{t\to\infty}N_0(x_0,t)=0$. Indeed, at the point $x_0$ the equation (2.4) takes the form
$$
  \frac{\partial N(x_0,t)}{\partial t}=-v_e(x_0)N(x_0,t) \quad \mbox{for}\;\; t>0,
$$
  so that $N(x_0,t)=N_0(x_0)e^{-tv_e(x_0)}$ for all $t\geq 0$, which implies that $\displaystyle
  \lim_{t\to\infty}N_0(x_0,t)=0$. This means that if the population does not migrate into a
  location but keeps migrating out of that location, then the population at that location will
  finally vanish.

  (6)\ \ If $v_e(x_0)=0$ and $v_i(x_0)>0$ for some $x_0\in\Omega$, then $N_0(x_0,t)$ is strictly
  monotone increasing in $t$. Indeed, at the point $x_0$ the equation (2.4) takes the form
$$
  \frac{\partial N(x_0,t)}{\partial t}=\int_{O}v(x_0,y)N(y,t)\rmd y>0
  \quad \mbox{for}\;\; t>0,
$$
  so that $N(x_0,t)$ is strictly monotone increasing. This means that if the population
  does not migrate out of a location but keeps constantly migrating to that location,
  then the population at that location keeps increasing.

\subsection{Solutions of the equation (2.3) in some simple situations}

\hskip 2em
  Our main interest is on asymptotic behavior of the solution $N(x,t)$ of the problem (2.3)--(2.4)
  as time $t$ goes to infinity. In this regards, the steady-state form of the equation (2.3), i.e.
  the equation
\begin{equation}
  \int_{\Omega}v(x,y)N^*(y)\rmd y=v_e(x)N^*(x), \quad x\in\Omega
\end{equation}
  must also be taken into account. As usual we call solution of this equation as {\em steady-state
  solution} of the problem (2.3)--(2.4).

  Before generally studying existence and uniqueness of a solution to the equation (2.6) and
  asymptotic behavior of the solution $N(x,t)$ of the problem (2.3)--(2.4), we first consider
  several special situations.

  (1)\ \ First consider {\em homogeneous migration}, i.e., the case where $v(x,y)=c=const.$ In this
  case, a simple computation shows that the solution of the equation (2.6) is given by $N^*(x)=
  C/|\Omega|$, where $C$ is an arbitrary constant, and the solution of the problem (2.3)--(2.4) is
  given by
$$
  N(x,t)=N_0(x)e^{-ct|\Omega|}+\frac{M_0}{|\Omega|}\Big(1-e^{-ct|\Omega|}\Big),
$$
  so that
\begin{equation}
  \lim_{t\to\infty}N(x,t)=\frac{M_0}{|\Omega|} \quad
  \mbox{uniformly for}\;\; x\in\Omega.
\end{equation}
  That is, if migration rate is homogeneous everywhere, then the population will finally
  homogeneously distributed in the habitat domain, regardless whether initially the population is
  homogeneously distributed or not.

  (2)\ \ Next consider balanced migration, i.e., the case where $v_i(x)=v_e(x)$ for all $x\in
  \Omega$. In this case, we can also easily verify that the solution of the problem (2.6) is given
  by $N^*(x)=C/|\Omega|$, where $C$ is an arbitrary constant. Later we shall see that the relation
  (2.7) is still valid for this case. Hence, for balanced migration, the population will finally
  homogeneously distributed in the habitat domain, regardless whether initially the
  population is homogeneously distributed or not.

  (3)\ \ We now consider the case that $v(x,y)$ depends only on $x$ and is independent of
  $y$: $v(x,y)=v(x)$, $\forall x,y\in\Omega$. In this case, the solution of the problem (2.7)
  is given by $N^*(x)=Cv(x)$, where $C$ is an arbitrary constant, and the solution of the problem
  (2.3)--(2.4) is given by
$$
  N(x,t)=N_0(x)e^{-at}+a^{-1}M_0v(x)\Big(1-e^{-at}\Big),
$$
  where $a=\displaystyle\int_{\Omega}v(x)\rmd x$, so that
\begin{equation*}
  \lim_{t\to\infty}N(x,t)=a^{-1}M_0v(x) \quad
  \mbox{uniformly for}\;\; x\in\Omega.
\end{equation*}
  Hence, in the case that the migration rate is determined only by the destination
  location and does not depend on the location of departure, the final distribution of
  the population is proportional to the migration rate function.

  (4)\ \ Consider another interesting case: $\Omega=\Omega_1\cup\Omega_2$, where
  $\Omega_1\cap\Omega_2=\emptyset$, and
$$
  v(x,y)=\left\{
\begin{array}{l}
  v_1(x), \quad \mbox{if}\;\; y\in\Omega_1,\\
  v_2(x), \quad \mbox{if}\;\; y\in\Omega_2,
\end{array}
\right.
  \quad x\in\Omega,
$$
  where $v_1$ and $v_2$ are given nonnegative functions in $\Omega$. We assume that $v$ is ergodic,
  which is equivalent to the following condition:
$$
  v_1(x)>0 \;\; \mbox{for}\;\; x\in\Omega_2 \quad \mbox{and} \quad
  v_2(x)>0 \;\; \mbox{for}\;\; x\in\Omega_1.
$$
  In this case, we have $v_e(x)=b_1\chi_{\Omega_1}(x)+b_2\chi_{\Omega_2}(x)$, where $\chi_S$
  denotes the characteristic function of the set $S\subseteq\Omega$, and $b_j=\displaystyle
  \int_{\Omega}v_j(x)\rmd x$, $j=1,2$. Let $a_{jk}=\displaystyle\int_{\Omega_j}v_k(x)\rmd x$, $j,k=1,2$.
  Note that the ergodicity assumption of $v$ implies that $a_{12}>0$ and $a_{21}>0$. We can easily
  verify that the solution of the equation (2.6) is given by
$$
  N^*(x)=\frac{a_{12}v_1(x)+a_{21}v_2(x)}{(a_{12}+a_{21})
  [b_1\chi_{\Omega_1}(x)+b_2\chi_{\Omega_2}(x)]}
$$
  (neglecting an arbitrary constant factor). By denoting $y_j(t)=\displaystyle\int_{\Omega_j}N(x,t)
  \rmd x$, $j=1,2$, we have
\begin{equation}
  N(x,t)=N_0(x)e^{-tv_e(x)}+v_1(x)\int_0^ty_1(s)e^{-(t-s)v_e(x)}ds
  +v_2(x)\int_0^ty_2(s)e^{-(t-s)v_e(x)}ds,
\end{equation}
  and $y_1,y_2$ satisfy the following system of differential equations:
$$
  \left\{
\begin{array}{l}
  y_1'(t)=-a_{21}y_1(t)+a_{12}y_2(t),
\\
  y_2'(t)=a_{21}y_1(t)-a_{12}y_2(t).
\end{array}
\right.
$$
  The solution of the above system is given by
\begin{equation}
  \left\{
\begin{array}{l}
  y_1(t)=\displaystyle y_1(0)e^{-(a_{12}+a_{21})t}+\frac{a_{12}M_0}{a_{12}+a_{21}}\Big(1-e^{-(a_{12}+a_{21})t}\Big),
\\ [0.3cm]
  y_2(t)=\displaystyle y_2(0)e^{-(a_{12}+a_{21})t}+\frac{a_{21}M_0}{a_{12}+a_{21}}\Big(1-e^{-(a_{12}+a_{21})t}\Big).
\end{array}
\right.
\end{equation}
  From (2.8) and (2.9) we can easily prove that
$$
  \lim_{t\to\infty}N(x,t)=M_0N^*(x)=\frac{M_0[a_{12}v_1(x)+a_{21}v_2(x)]}
  {(a_{12}+a_{21})[b_1\chi_{\Omega_1}(x)+b_2\chi_{\Omega_2}(x)]}
$$
  uniformly for $x\in\Omega$.
\medskip

  (5)\ \ The above example can be extended to the general $n$-piece case. We make a short
  discussion to this case here. Thus we assume that $\Omega=\displaystyle\bigcup_{k=1}^n
  \Omega_k$, where $\Omega_j\cap\Omega_k=\emptyset$ for any pair $1\leq j,k\leq n$ with
  $j\neq k$, and
$$
  v(x,y)=v_k(x) \quad \mbox{for}\;\;  x\in\Omega,\;\; y\in\Omega_k, \;\; k=1,2,\cdots,n.
$$
  We call such $v$ a {\em piecewise semi-constant migration function}. Let
$$
  a_{jk}=\int_{\Omega_j}\!\!v_k(x)\rmd x, \quad b_k=\sum_{j=1}^n a_{jk}=\int_{\Omega}v_k(x)\rmd x,
  \quad j,k=1,2,\cdots,n,
$$
  and
$$
  y_j(t)=\int_{\Omega_j}\!\!N(x,t)\rmd x,  \quad j=1,2,\cdots,n.
$$
  Clearly, $v_e(x)=\displaystyle\sum_{k=1}^nb_k\chi_{\Omega_k}(x)$, and $y_1,y_2,\cdots,y_n$ satisfy
  the following system of ordinary differential equations:
\begin{equation}
  y_j'(t)=\sum_{k=1}^na_{jk}y_k(t)-b_jy_j(t), \quad j=1,2,\cdots,n.
\end{equation}
  Notice that when this system of equations is solved, the equation (2.3) reduces into the
  following simpler equation:
$$
  \frac{\partial N(x,t)}{\partial t}=\sum_{k=1}^ny_k(t)v_k(x)-v_e(x)N(x,t).
$$
  Hence the solution of the problem (2.3)--(2.4) is given by
\begin{equation}
  N(x,t)=N_0(x)e^{-tv_e(x)}+\sum_{k=1}^nv_k(x)\int_0^ty_k(s)e^{-(t-s)v_e(x)}ds.
\end{equation}
  To study asymptotic behavior of solutions of the system (2.10), let $A=(a_{jk})_{n\times n}$ and
  $B=\mbox{\rm diag}(b_1,b_2,\cdots,b_n)$, where $b_k=\displaystyle\sum_{j=1}^n a_{jk}$, $k=1,2,
  \cdots,n$. As before we assume that $v$ is ergodic. Then it is easy to see that the matrix $A$ is
  irreducible, i.e., for any pair of distict indices $1\leq j,k\leq n$, there exists a corresponding
  chain of indices $k=j_0,j_1,\cdots,j_l=j$ such that $a_{j_sj_{s-1}}>0$ for all $s=1,2,\cdots,l$.
  By the famous Perron-Frobenious theorem, it follows that $0$ is an algebraically simple eigenvalue
  of $A-B$ with a strictly positive eigenvector, and there exists constant $\kappa>0$ such that
  for any nonzero eigenvalue $\lambda$ of $A-B$ we have $\mbox{\rm Re}\lambda\leq-\kappa$. Let
  $\xi=(\xi_1,\xi_2,\cdots,\xi_n)\in\mathbb{R}^n$ be the strictly positive eigenvector of $A-B$
  corresponding to the $0$ eigenvalue, i.e., $A\xi^T=B\xi^T$, normalized such that $\displaystyle
  \sum_{k=1}^n \xi_k=1$. Let $c=\displaystyle\sum_{j=1}^n y_j(0)$ and $z_j(t)=y_j(t)-c\xi_j$,
  $j=1,2,\cdots,n$, $t\geq 0$. It follows that $(z_1(t),z_2(t),\cdots,z_n(t))$ is a solution of
  the system of equations (2.10) satisfying
\begin{equation}
  \sum_{j=1}^n z_j(t)=0 \quad \mbox{for}\;\; t\geq 0.
\end{equation}
  Moreover, for any $0<\kappa'<\kappa$ there exists corresponding constant $C>0$ such that
\begin{equation}
  \sum_{j=1}^n |z_j(t)|\leq Ce^{-\kappa' t}\sum_{j=1}^n |y_j(0)-c\xi_j|
  \quad \mbox{for}\;\; t\geq 0.
\end{equation}
  In other words, the solution of (2.10) has the following expression:
\begin{equation}
  y_j(t)=c\xi_j+z_j(t), \quad t\geq 0,\;\; j=1,2,\cdots,n,
\end{equation}
  where $c=\displaystyle\sum_{j=1}^n y_j(0)$ and $(z_1(t),z_2(t),\cdots,z_n(t))$ is a solution of
  the system of equations (2.10) satisfying (2.12) and (2.13). Substituting the expression (2.14)
  into (2.11), we get the following result:
\medskip

  {\bf Theorem 2.3}\ \ {\em Let $v(x,y)$ be a piecewise semi-constant migration function and let
  $N=N(x,t)$ be the solution of the problem $(2.3)$--$(2.4)$. Assume that $v$ is ergodic. Then
\begin{equation*}
  N(x,t)=M_0N^*(x)+[N_0(x)-M_0N^*(x)]e^{-tv_e(x)}+h_0(x,t)
  \quad \mbox{for}\;\; x\in\Omega,\;\, t\geq 0,
\end{equation*}
  where $N^*(x)=\displaystyle v_e^{-1}(x)\sum_{k=1}^n\xi_kv_k(x)$, $M_0=\displaystyle\int_{\Omega}
  N_0(x)\rmd x$, and $h_0$ is a function uniquely determined by $N_0$ and possesses the following property:
$$
  \int_{\Omega}|h_0(x,t)|\rmd x\leq Ce^{-\kappa' t}\int_{\Omega}|N_0(x)-M_0N^*(x)|\rmd x
  \quad \mbox{for} \;\; t\geq 0,
$$
  where $0<\kappa'<\kappa$ and $C$ is a positive constant depending on $\kappa'$. $\quad\Box$}

\subsection{Asymptotic behavior of solutions in the ergodic migration case}

\hskip 2em
  We now turn to consider the case where the migration rate function $v$ is continuous and ergodic.
  We first consider the steady-state equation (2.6). For it we have the following result:
\medskip

  {\bf Theorem 2.4}\ \ {\em Assume that $v\in C(\overline{\Omega}\times\overline{\Omega})$
  and is ergodic. The equation $(2.6)$ has a unique solution $N^*\in C(\overline{\Omega})$
  satisfying the normalization condition $\displaystyle\int_{\Omega}N^*(x)\rmd x=1$. This solution
  is strictly positive in $\overline{\Omega}$, i.e., $N^*(x)>0$ for all $x\in\overline{\Omega}$,
  and all the other solutions of $(2.6)$ are all proportional to this solution.}
\medskip

  {\em Proof}:\ \ Let $a(x,y)=v(x,y)/v_e(y)$, $x,y\in\Omega$ and $\varphi(x)=
  v_e(x)N^*(x)$, $x\in\Omega$. Then the equation (2.6) is transformed into the following equivalent
  equation:
\begin{equation}
  \int_{\Omega}a(x,y)\varphi(y)\rmd y=\varphi(x), \quad x\in\Omega.
\end{equation}
  Consider the operator $K:\varphi\mapsto\displaystyle\int_{\Omega}a(x,y)\varphi(y)\rmd y$, $\forall
  \varphi\in C(\overline{\Omega})$. It is clear that $K$ is a positive compact linear operator in
  $C(\overline{\Omega})$. We claim that the spectral radius of $K$ is equal to $1$: $r(K)=1$. To
  see this, we first note that since $\displaystyle\int_{\Omega}a(x,y)\rmd x=1$ for all $y\in
  \overline{\Omega}$, we have $K'1_{\Omega}=1_{\Omega}$, where $K'$ denotes the conjugate
  of $K$ and $1_{\Omega}$ denotes the unit-valued constant function in $\overline{\Omega}$.
  This shows that $1$ is an eigenvalue of $K'$, so that $1\in\sigma(K')=\sigma(K)$, which implies
  that $r(K)\geq 1$. To prove the opposite inequality we let $\lambda\in\mathbb{C}$ be such that
  $|\lambda|>1$. Since for any $\varphi\in L^1(\Omega)$,
$$
  \|K\varphi\|_{L^1}=\int_{\Omega}\!\Big|\!\int_{\Omega}\!a(x,y)\varphi(y)\rmd y\Big|\rmd x
  \leq\int_{\Omega}\!\Big(\!\int_{\Omega}\!a(x,y)\rmd x\Big)|\varphi(y)|\rmd y=\|\varphi\|_{L^1},
$$
  we see that $\|K\|_{\mathscr{L}(L^1(\Omega))}\leq 1$, so that for any $|\lambda|>1$, as a
  bounded linear operator in $L^1(\Omega)$, $\lambda I-K$ has bounded inverse
$$
  (\lambda I-K)^{-1}=\lambda^{-1}(I-\lambda^{-1}K)^{-1}
  =\lambda^{-1}\sum_{n=0}^{\infty}(\lambda^{-1}K)^n,
$$
  where convergence of the series in $\mathscr{L}(L^1(\Omega))$ follows from the fact that
  $\|\lambda^{-1}K\|_{\mathscr{L}(L^1(\Omega))}\leq|\lambda|^{-1}<1$. Moreover, the above
  equalities also show that
$$
  (\lambda I-K)^{-1}=\lambda^{-1}I+\lambda^{-1}K(\lambda I-K)^{-1}.
$$
  Since $(\lambda I-K)^{-1}$ is a bounded linear operator in $L^1(\Omega)$ so that it maps
  $C(\overline{\Omega})$ to $L^1(\Omega)$ boundedly and clearly $K$ maps $L^1(\Omega)$ to
  $C(\overline{\Omega})$ boundedly, from the above equality we see that the restriction of
  $(\lambda I-K)^{-1}$ in $C(\overline{\Omega})$ is a bounded linear operator in this space. This
  shows that $\lambda\in\rho(K)$ whenever $|\lambda|>1$, so that $r(K)\leq 1$. Hence our claim is
  true. Now, since $K$ is irreducible (by Lemma 2.2), it follows by a generalization of the
  well-known Krein-Rutman theorem, or more precisely Theorem 5.2 in Chapter V of \cite{Sch}, that
  the spectral radius $r(K)=1$ is a simple eigenvalue of $K$ corresponding to a strictly positive
  eigenvector. Let $\varphi\in C(\overline{\Omega})$ be a strictly positive eigenvector of $K$
  corresponding to the eigenvalue $r(K)=1$. Then $\varphi$ is a strictly positive solution of the
  equation (2.15). Moreover, simplicity of $1$ as an eigenvalue of $K$ implies that all other
  solutions of (2.15) are constant multiples of $\varphi$. Now set $N^*(x)=c\varphi(x)/v_e(x)$,
  $x\in\Omega$, where $c$ is the normalization factor. Then $N^*$ is a strictly positive solution
  of the equation (2.6). Moreover, from the above assertion on the equation (2.15) we easily see
  that any other solution of (2.6) is a multiple of $N^*$ by a constant factor. This proves Theorem
  2.4. $\quad\Box$
\medskip

  Later on, we shall always use the notation $N^*$ to denote the unique normalized solution of
  (2.6) ensured by the above theorem, i.e., $N^*$ is the unique solution of the following problem:
\begin{equation}
\left\{
\begin{array}{l}
  \displaystyle\int_{\Omega}v(x,y)N^*(y)\rmd y=v_e(x)N^*(x), \quad x\in\Omega,
\\ [0.3cm]
  \displaystyle\int_{\Omega}N^*(x)\rmd x=1.
\end{array}
\right.
\end{equation}
  By Theorem 2.1 we know that $N^*(x)>0$ for all $x\in\overline{\Omega}$, and for any other
  solution $\tilde{N}^*$ of the equation (2.6) there exists a corresponding constant $c$ such that
  $\tilde{N}^*(x)=cN^*(x)$ for all $x\in\overline{\Omega}$.
\medskip

  For the solution of the problem (2.3)--(2.4), we have the following result:
\medskip

  {\bf Theorem 2.5}\ \ {\em Assume that $v\in C(\overline{\Omega}\times\overline{\Omega})$ and is
  ergodic. Let $N^*$ be as in Theorem 2.4. For $N_0\in C(\overline{\Omega})$, let $M_0=
  \displaystyle\int_{\Omega}N_0(x)\rmd x$ and $N=N(x,t)$ be the solution of the problem
  $(2.3)$--$(2.4)$. Then
\begin{equation}
  N(x,t)=M_0N^*(x)+h(x,t) \quad \mbox{for}\;\; x\in\Omega,\;\, t\geq 0,
\end{equation}
  where $h$ is a function uniquely determined by $N_0$ possessing the following properties:

  $(1)$\ \ $\displaystyle\int_{\Omega}h(x,t)\rmd x=0$ for all $t\geq 0$,

  $(2)$\ \ $\displaystyle\sup_{x\in\Omega}|h(x,t)|\leq Ce^{-\kappa t}\sup_{x\in\Omega}
  |N_0(x)-M_0N^*(x)|$ for all $t\geq 0$, and

  $(3)$\ \ $\displaystyle\int_{\Omega}|h(x,t)|\rmd x\leq Ce^{-\kappa t}\int_{\Omega}
  |N_0(x)-M_0N^*(x)|\rmd x$ for all $t\geq 0$.

\noindent
  Here $\kappa$ and $C$ are positive constants depending only on the function $v$.}
\medskip

  {\em Proof}:\ \ We denote
$$
  X=C(\overline{\Omega}), \quad X_0=\mbox{\rm span}\{N^*\}, \quad
  X_1=\Big\{\varphi\in X:\int_{\Omega}\varphi(x)\rmd x=0\Big\}.
$$
  It is clear that both $X_0$ and $X_1$ are closed linear subspaces of $X$, and $X=X_0\oplus X_1$.
  Let $H=A-B:X\to X$, where $A:X\to X$ and $B:X\to X$ are respectively defined as follows:
\begin{equation*}
  A\varphi(x)=\int_{\Omega}v(x,y)\varphi(y)\rmd y, \quad
  B\varphi(x)=v_e(x)\varphi(x), \quad \forall\varphi\in X.
\end{equation*}
  Theorem 2.4 ensures that $X_0=\ker(H)$. Moreover, a simple computation shows that $\mbox{rg}(H)
  \subseteq X_1$. Hence, both $X_0$ and $X_1$
  are invariant subspaces of $H$. It follows that
\begin{equation}
  \sigma(H)=\sigma(H|_{X_0})\cup\sigma(H|_{X_1})=\{0\}\cup\sigma(H|_{X_1}),
\end{equation}
  and
\begin{equation}
  e^{tH}=e^{tH|_{X_0}}\oplus e^{tH|_{X_1}}=id_{X_0}\oplus e^{tH|_{X_1}} \quad \forall t\geq 0,
\end{equation}
  i.e., for any $\varphi\in X$, letting $\varphi=\varphi_0+\varphi_1$, where $\varphi_i\in X_i$,
  $i=0,1$, we have $e^{tH}\varphi=\varphi_0+e^{tH|_{X_1}}\varphi_1$ for all $t\geq 0$. Now, for
  any $N_0\in X=C(\overline{\Omega})$, the decomposition $X=X_0\oplus X_1$ implies that
$$
  N_0=M_0N^*+h_0,
$$
  where $M_0=\displaystyle\int_{\Omega}N_0(x)\rmd x$ and $h_0=N_0-M_0N^*\in X_1$. Hence, by letting
  $h(\cdot,t)=e^{tH|_{X_1}}h_0$, we have $h(\cdot,t)\in X_1$ for all $t\geq 0$ and, by (2.19),
  we have
$$
  N(\cdot,t)=e^{tH}N_0=M_0N^*+h(\cdot,t) \quad \forall t\geq 0.
$$
  This proves (2.17) as well as the assertion (1).

  Next we prove that $s(H)=0$. Since $0$ is an eigenvalue of $H$ so that $0\in\sigma(H)$, it is
  clear that $s(H)\geq 0$. Hence we only need to prove that also $s(H)\leq 0$. To this end, by
  Lemma 1.9 in Chapter VI of \cite{EN}, it suffices to prove that for any $\lambda>0$ we have
  $\lambda\in\rho(H)$ and $R(\lambda,H)\geq 0$. This follows from a similar argument as in the
  proof of Theorem 2.4. Indeed, for any $\lambda>0$ and $f\in X$, the equation
$$
  (\lambda I-H)u=f \quad \mbox{or} \quad  (\lambda I+B)u-Au=f
$$
  can be explicitly rewritten as follows:
$$
  [\lambda+v_e(x)]u(x)-\int_{\Omega}v(x,y)u(y)\rmd y=f(x), \quad x\in\Omega.
$$
  Let $\hat{v}(x,y)=v(x,y)/[\lambda+v_e(y)]$ and $\hat{u}(x)=[\lambda+v_e(x)]u(x)$. Then the
  above equation can be rewritten as follows:
\begin{equation}
  \hat{u}(x)-\int_{\Omega}\hat{v}(x,y)\hat{u}(y)\rmd y=f(x), \quad x\in\Omega.
\end{equation}
  Let $\hat{K}$ be the following bounded linear operator from $L^1(\Omega)$ to $X=
  C(\overline{\Omega})$:
$$
  \varphi\mapsto \hat{K}\varphi=\int_{\Omega}\hat{v}(x,y)\varphi(y)\rmd y
  \quad \mbox{for}\;\; \varphi\in L^1(\Omega),
$$
  and let $\theta=\displaystyle\max_{y\in\overline{\Omega}}v_e(y)/[\lambda+v_e(y)]$. Clearly
  $0<\theta<1$ and $\displaystyle\|\hat{K}\|_{\mathscr{L}(L^1(\Omega))}\leq\theta< 1$. Hence a
  similar argument as in the proof of Theorem 2.4 shows that for any $f\in X$ the equation (2.20)
  has a unique solution $\hat{u}\in X$, given by $\hat{u}=Lf$, where $L=\displaystyle
  \sum_{n=0}^{\infty}\hat{K}^n=I+\hat{K}\sum_{n=0}^{\infty}\hat{K}^n$ is a bounded linear operator
  in $X$. Since the mapping $\hat{u}\mapsto u=\hat{u}/(\lambda+v_e)$ is clearly a bounded linear
  operator in $X$, it follows that for any $\lambda>0$ we have $\lambda\in\rho(H)$. Moreover, from
  the above expression we easy see that if $f\geq 0$ then $\hat{u}\geq 0$ so that also $u\geq 0$,
  i.e., $R(\lambda,H)=(\lambda I-H)^{-1}\geq 0$ for any $\lambda>0$. This proves the desired
  assertion and, consequently, we have $s(H)=0$.

  In what follows we prove that $s(H_{X_1})<0$. First we note that (2.18) implies that
  $\sigma(H_{X_1})\subseteq\sigma(H)$. Moreover, since $0$ is a simple eigenvalue of $H$, we have
  $0\not\in\sigma(H_{X_1})$. Hence $\sigma(H_{X_1})=\sigma(H)\backslash\{0\}$. In the sequel we
  prove that $(\sigma(H)\backslash\{0\})\cap i\mathbb{R}=\emptyset$.

  We first prove that if $\lambda\in\partial\sigma(H)\backslash\mathbb{R}$ then $\lambda$ is an
  eigenvalue of $H$. Let $\lambda\in\partial\sigma(H)$ and $\lambda=\mu+i\nu$, $\nu\neq 0$. Since
  $\lambda\in\partial\sigma(H)$, by a well-known result we see that $\lambda$ is an approximate
  eigenvalue of $H$, i.e., there exists a sequence $\{u_n\}_{n=1}^{\infty}\subseteq X$ with
  $\|u_n\|_X=1$ ($n=1,2,\cdots$), such that $\displaystyle\lim_{n\to\infty}\|\lambda u_n-Hu_n\|_X
  =0$ (cf. Proposition 1.10 in Chapter IV of \cite{EN}). Let $f_n=\lambda u_n-Hu_n$ ($n=1,2,\cdots$).
  Then $\displaystyle\lim_{n\to\infty}\|f_n\|_X=0$ and for each $n\in\mathbb{N}$, $u_n$ is a
  solution of the following integral equation:
$$
  [\lambda+v_e(x)]u_n(x)-\int_{\Omega}v(x,y)u_n(y)\rmd y=f_n(x), \quad x\in\Omega.
$$
  Let $c=\displaystyle\min_{x\in\overline{\Omega}}|\lambda+v_e(x)|$. Since $\mbox{\rm Im}\lambda
  \neq 0$ and $v_e$ is a real-valued continuous function in $\overline{\Omega}$, we have $c>0$.
  Let $\hat{v}(x,y)=v(x,y)/[\lambda+v_e(y)]$, $\hat{u}_n=(\lambda+v_e)u_n/\|(\lambda+v_e)u_n\|_X$
  and $\hat{f}_n=f_n/\|(\lambda+v_e)u_n\|_X$. Then the above equation can be rewritten as follows:
\begin{equation*}
  \hat{u}_n(x)-\int_{\Omega}\hat{v}(x,y)\hat{u}_n(y)\rmd y=\hat{f}_n(x), \quad x\in\Omega.
\end{equation*}
  Note that $\|\hat{u}_n\|_X=1$ and $\|\hat{f}_n\|_X\leq c^{-1}\|f_n\|_X$, $n=1,2,\cdots$, so that
  $\displaystyle\lim_{n\to\infty}\|\hat{f}_n\|_X=0$. Hence, the above equations imply that $1$ is
  an approximate eigenvalue of the compact linear operator $K:X\to X$ defined as before. By the
  Riesz-Schauder theory for compact linear operators in Banach spaces, it follows that $1$ is an
  eigenvalue of $K$ and, consequently, $\lambda$ is an eigenvalue of $H$. This proves the desired
  assertion.

  As an immediate consequence of the above assertion it follows that
$$
  \sigma_b(H):=\{\lambda\in\sigma(H):\mbox{\rm Re}\lambda=s(H)\}=P\sigma(H)\cap i\mathbb{R},
$$
  where $P\sigma(H)$ denotes the point spectrum of $H$. By Lemma 2.2, ergodicity of the kernel $v$
  ensures that the semigroup $e^{tH}$ is irreducible. Hence, by Theorem 3.6 in Chapter B and
  Theorem 3.8 in C of \cite{Nag}, it follows that $P\sigma(H)\cap i\mathbb{R}$ is cyclic, i.e.,
  $P\sigma(H)\cap i\mathbb{R}$ is a additive subgroup of $i\mathbb{R}$. Since $H$ is a bounded
  linear operator so that its spectrum $\sigma(H)$ is a bounded subset of $\mathbb{C}$, this implies
  that $\sigma_b(H)=P\sigma(H)\cap i\mathbb{R}=\{0\}$. Hence $(\sigma(H)\backslash\{0\})\cap
  i\mathbb{R}=\emptyset$, as desired.

  As a consequence of the above assertion, we have $\sigma(H_{X_1})\cap i\mathbb{R}=\emptyset$.
  Since $\sigma(H_{X_1})\subseteq\sigma(H)\subseteq\{\lambda\in\mathbb{C}:\mbox{\rm Re}\lambda\leq
  0\}$ and it is a bounded closed subset of $\mathbb{C}$, it follows that $s(H_{X_1})=
  \sup\{\mbox{\rm Re}\lambda:\lambda\in\sigma(H_{X_1})\}<0$. Hence, by arbitrarily choosing a
  number $\kappa<0$ such that $s(H_{X_1})<\kappa<0$, we obtain the assertion (2).

  Finally, to prove the assertion (3) we let $Y=L^1(\Omega)$, $Y_0=X_0$, and
$$
  Y_1=\Big\{\varphi\in Y:\int_{\Omega}\varphi(x)\rmd x=0\Big\}.
$$
  Clearly, $Y=Y_0\oplus Y_1$. By regarding $H$ as a bounded linear operator in $Y$, we have
\begin{equation}
  \sigma(H)=\sigma(H|_{Y_0})\cup\sigma(H|_{Y_1})=\{0\}\cup\sigma(H|_{Y_1}),
\end{equation}
  and
\begin{equation}
  e^{tH}=e^{tH|_{Y_0}}\oplus e^{tH|_{Y_1}}=id_{Y_0}\oplus e^{tH|_{Y_1}}.
\end{equation}
  Since $v_e(x)$ and $v(x,y)$ are positive continuous functions, it can be easily seen that
  $\ker H\subseteq C(\overline{\Omega})$ and, consequently, $0$ is still a simple eigenvalue of
  $H$ when regarding it as a bounded linear operator in $Y$. This implies that $0\not\in
  \sigma(H|_{Y_1})$. Using this fact and a similar argument as above we can still prove that
  $s(H|_{Y_1})<0$, which yields the assertion (3). This completes the proof of Theorem 2.5.
  $\quad\Box$
\medskip

  Note that Theorem 2.5 implies that if $N_0\in C(\overline{\Omega})$ then
\begin{equation}
  \lim_{t\to\infty}\sup_{x\in\Omega}|N(x,t)-M_0N^*(x)|=0,
\end{equation}
  and if $N_0\in L^1(\Omega)$ then
\begin{equation}
  \lim_{t\to\infty}\int_{\Omega}|N(x,t)-M_0N^*(x)|\rmd x=0.
\end{equation}
  Moreover, the convergence is exponentially fast.

  Hence, Theorems 2.4 and 2.5 show that in the proliferation-stationary case, if the population
  system is ergodic then in addition to the property that total amount of the population keeps
  constant (which is a consequence of the stationary proliferation of the system), the population
  of the system will finally distributed in the whole habitat domain in an intrisic way uniquely
  determined by the migration function. More precisely, the final distribution of the population
  is proportional to a specific function $N^*$ which is everywhere positive in the habitat domain
  $\Omega$. This function $N^*$ depends only on the migration function $v$ of the system and is
  independent of the specific initial distribution of the population $N_0$. Moreover, the proportion
  coefficient is uniquely determined by the initial total amount $M_0=\displaystyle\int_{\Omega}
  N_0(x)\rmd x$ of the population and is also independent of the specific initial distribution $N_0$
  of the population.

\subsection{Proliferation non-stationary population migration model}

\hskip 2em
  In this subsection we extend the result of the above subsection to the proliferation
  non-stationary case. Let $r=r(x)$, $x\in\Omega$, be the proliferation rate function, i.e., for
  every $x\in\Omega$, $r(x)$ is the proliferation rate (=birth rate minus death rate) of the
  population at location $x$. Assume that this function is not identically vanishing in $\Omega$.
  Then the equation (2.4) should be replaced by the following equation:
\begin{equation}
  \frac{\partial N(x,t)}{\partial t}=\int_{\Omega}v(x,y)N(y,t)\rmd y-v_e(x)N(x,t)+
  r(x)N(x,t), \quad  x\in\Omega,\;\; t>0.
\end{equation}
  To study this equation, we first make a short discussion to the spectral of a general class of
  bounded linear operators in $C(\overline{\Omega})$.

  Let $a\in C(\overline{\Omega}\times\overline{\Omega})$ be a given function such that $a(x,y)\geq
  0$ for all $x,y\in\overline{\Omega}$, and let $b\in C(\overline{\Omega})$. We denote by $A$ the
  following bounded linear operator in $C(\overline{\Omega})$:
\begin{equation}
  A\varphi(x)=\int_{\Omega}a(x,y)\varphi(y)\rmd y+b(x)\varphi(x),
  \quad \mbox{for}\;\; \varphi\in C(\overline{\Omega}).
\end{equation}
  We use the notation $s(A)$ to denote the {\em spectral bound} of $A$, i.e., $s(A)=\displaystyle
  \max_{\lambda\in\sigma(A)}\mbox{\rm Re}\lambda$.
\medskip

  {\bf Lemma 2.6}\ \ {\em Let $\Omega$ be a bounded domain in $\mathbb{R}^n$ and $A\in
  \mathscr{L}(C(\overline{\Omega}))$ be as above. We have the following assertions:

  $(1)$ $s(A)\in\sigma(A)$.

  $(2)$ $\{b(x):x\in\overline{\Omega}\}\subseteq\sigma(A)$.

  $(3)$ $\sigma(A)\backslash\{b(x):x\in\overline{\Omega}\}$ contains only poles of finite algebraic
  multiplicity of the resolvent $R(\cdot,A)$ so that they are eigenvalues of $A$ of finite algebraic
  multiplicity.

  $(4)$ The following estimates hold:
\begin{equation}
  \max_{x\in\overline{\Omega}}b(x)\leq s(A)\leq\max_{y\in\overline{\Omega}}\Big[\int_{\Omega}
  a(x,y)\rmd x+b(y)\Big].
\end{equation}

  $(5)$ If $\lambda>s(A)$ then $R(\lambda,A)=(\lambda I-A)^{-1}$ is a positive operator.

  $(6)$ There exists a positive probability measure $\rmd\mu\in\mathcal{M}(\Omega)$ such that
  $A'\rmd\mu=s(A)\rmd\mu$, where $\mathcal{M}(\Omega)$ is the Banach space of finite Borel measures
  in $\Omega$, and $A'$ denotes the transpose of $A$.

  $(7)$ Assume that $a$ is ergodic  and $s(A)>\displaystyle\max_{x\in\overline{\Omega}}b(x)$. Then
  $s(A)$ is a simple eigenvalue of $A$ with a strictly positive eigenvector. Moreover,  $s(A)$ is
  a dominant spectral value of $A$, i.e., $\mbox{\rm Re}\lambda<s(A)$ for any $\lambda\in\sigma(A)
  \backslash\{s(A)\}$.}
\medskip

  {\em Proof}:\ \ It is easy to see that $A$ generates a positive uniformly continuous semigroup in
  $C(\overline{\Omega})$. Hence the assertion (1) follows from a well-known result for positive
  $C_0$-semigroups, cf. Theorem 1.10 in Chapter VI of \cite{EN}. To prove the assertion (2) we
  assume that $x_0\in\overline{\Omega}$ and let $\lambda=b(x_0)$. Since $b$ is continuous at $x_0$, for
  each $n\in\mathbb{N}$ there exists a
  corresponding $\delta_n>0$ such that $|b(x)-b(x_0)|<1/n$ for all $x\in B(x_0,\delta_n)\cap
  \overline{\Omega}$. Moreover, we can assume that $\delta_n\to 0$ as $n\to\infty$. For every
  $n\in\mathbb{N}$, let $\varphi_n\in C(\overline{\Omega})$ be a cutoff function such that
$$
  0\leq\varphi_n\leq 1, \quad \varphi_n(x_0)=1, \quad \mbox{and} \quad
  \varphi_n(x)=0 \;\; \mbox{for}\;\; x\in\overline{\Omega}\backslash B(x_0,\delta_n).
$$
  Clearly, $\|\varphi_n\|_{\infty}=1$, $n=1,2,\cdots$, and
$$
  \|(\lambda I-A)\varphi_n\|_{\infty}=\max_{x\in\overline{\Omega}}
  \Big|[b(x_0)-b(x)]\varphi_n(x)-\int_{\Omega}a(x,y)\varphi_n(y)\rmd y\Big|
  \leq\frac{1}{n}+C\delta_n^2\to 0 \quad \mbox{as}\;\; n\to\infty.
$$
  Hence $\lambda=b(x_0)$ is an approximate eigenvalue of $A$, so that $\lambda\in\sigma(A)$. This
  proves the assertion (2).

  To prove the assertion (3) we let $M_b$ be the multiplication operator related to the function
  $b$, i.e., $M_b\varphi=b\varphi$ for $\varphi\in C(\overline{\Omega})$, and let $K$ be the
  integral operator with kernel $a$, i.e., $K\varphi(x)=\displaystyle\int_{\Omega}a(x,y)\varphi(y)
  \rmd y$ for $\varphi\in C(\overline{\Omega})$. Then $A=M_b+K$. Since $K$ is a compact operator,
  we have $\sigma_{ess}(A)=\sigma_{ess}(M_b)=\{b(x):x\in\overline{\Omega}\}$. Hence $\sigma(A)
  \backslash\{b(x):x\in\overline{\Omega}\}=\sigma(A)\backslash\sigma_{ess}(A)$. Since $\sigma(A)
  \backslash\sigma_{ess}(A)$ contains only poles of $R(\cdot,A)$ of finite algebraic multiplicity
  (cf. Section 1.20 in Chapter IV of \cite{EN}), we obtain the assertion (3).

  Next, the first inequality in (2.27) is an immediate consequence of the assertion (2) (it also
  easily follows from Corollary 1.11 in Chapter VI and Proposition 4.2 in Chapter I of \cite{EN}).
  We now prove the second inequality in (2.27). Let $\lambda$ be a complex number with real part
  greater than the number on the right-hand side of the second inequality in (2.27). Then
$$
  \mbox{\rm Re}\lambda-b(y)>\int_{\Omega}a(x,y)\rmd x\geq 0
  \quad \mbox{for}\;\; y\in\overline{\Omega}
$$
  and $\alpha=\displaystyle\max_{y\in\overline{\Omega}}\Big(\int_{\Omega}a(x,y)\rmd x\Big/
  [\mbox{\rm Re}\lambda-b(y)]\Big)<1$. For any given $f\in C(\overline{\Omega})$, by making the
  unknown variable transformation $\psi(x)=[\lambda-b(x)]\varphi(x)$, the equation $(\lambda I-A)
  \varphi=f$ or
\begin{equation}
  [\lambda-b(x)]\varphi(x)-\int_{\Omega}a(x,y)\varphi(y)\rmd y=f(x)
\end{equation}
  is transformed into the following equivalent equation:
\begin{equation}
  \psi(x)-\int_{\Omega}\frac{a(x,y)}{\lambda-b(y)}\psi(y)\rmd y=f(x).
\end{equation}
  Let $K_{\lambda}$ be as before, i.e., it denotes the integral operator defined by the expression
  on the left-hand side of the above equation. It is easy to see that
  $\|K_{\lambda}\|_{\mathscr{L}(L^1(\Omega))}\leq\alpha<1$. Hence by a similar argument as in the
  proof of Theorem 2.4 we see that $I-K_{\lambda}$, as a bounded linear operator in
  $C(\overline{\Omega})$, has bounded inverse in this space. This means that for any $f\in
  C(\overline{\Omega})$ the equation (2.31) has a unique solution $\psi\in C(\overline{\Omega})$,
  and the solution operator $f\mapsto\psi$ is a bounded linear operator in $C(\overline{\Omega})$.
  Since the mapping $\psi\mapsto\varphi=\psi/(\lambda-b)$ is also a bounded
  linear operator in $C(\overline{\Omega})$, it follows that for any $f\in C(\overline{\Omega})$
  the equation (2.29) has a unique solution $\varphi\in C(\overline{\Omega})$, and the solution
  operator $f\mapsto\varphi$ is a bounded linear operator in $C(\overline{\Omega})$. Hence,
  $\lambda\in\rho(A)$. This proves the second inequality in (2.27). This proves the assertion (4).

  The assertion (5) follows from Lemma 1.9 in Chapter VI of \cite{EN}. The assertion (6) follows
  from Theorem 1.6 in Chapter B-III of \cite{Nag}. To prove the assertion (7) we let $M_b$ and
  $K$ be as before. Then the relation $A=M_b+K$ and compactness of $K$ imply that $\omega_{ess}(A)
  =\omega_{ess}(M_b)=\displaystyle\max_{x\in\overline{\Omega}}b(x)$ (cf. Proposition 2.12 in
  Chapter IV of \cite{EN}). Hence, the assumption $s(A)>\displaystyle\max_{x\in\overline{\Omega}}
  b(x)$ implies that $s(A)>\omega_{ess}(A)$, so that by Corollary 3.16 in Chapter C-III of
  \cite{Nag} we see that $s(A)$ is a dominant spectral value of $A$ (This also follows from
  Corollary 2.13 in Chapter C-III of \cite{Nag}). Moreover, by applying Theorem
  3.8 in Chapter C-III of \cite{Nag} to the operator $A_1=A-s(A)I$ (which satisfies the condition
  $s(A_1)=0$) we see that $A$ has a positive eigenvector corresponding to the eigenvalue $s(A)$.
  Moreover, ergodicity of the function $a$ implies that the semigroup generated by $A$ is
  irreducible, so that by applying Proposition 3.5 in Chapter C-III of \cite{Nag} we see that this
  positive eigenvector of $A$ is strictly positive. This proves the assertion (7). The proof of
  Lemma 2.6 is complete. $\quad\Box$
\medskip

  We now use Lemma 2.6 to study asymptotic stability of solutions of the equation (2.25). Let
  $H$ be the bounded linear operator in $C(\overline{\Omega})$ defined by
\begin{equation}
  H\varphi(x)=\int_{\Omega}v(x,y)\varphi(y)\rmd y-v_e(x)\varphi(x)+r(x)\varphi(x),
  \quad  x\in\Omega,\;\; \quad \mbox{for}\;\; \varphi\in C(\overline{\Omega}).
\end{equation}
  Let $s(H)$ be the spectral bound of $H$, and $H':\mathcal{M}(\Omega)\to\mathcal{M}(\Omega)$
  be the transpose of $H$.
\medskip

  {\bf Theorem 2.7}\ \ {\em We have the following assertions:

  $(1)$ $\displaystyle\max_{x\in\overline{\Omega}}[r(x)-v_e(x)]\leq s(H)\leq
  \max_{x\in\overline{\Omega}}r(x)$.

  $(2)$ There exists a positive probability measure $\rmd\mu\in\mathcal{M}(\Omega)$ such that
  $H'\rmd\mu=s(H)\rmd\mu$, where $\mathcal{M}(\Omega)$ is the Banach space of finite Borel measures
  in $\Omega$, and $H'$ denotes the transpose of $H$.

  $(3)$ Assume that $v$ is ergodic and $s(H)>\displaystyle\max_{x\in\overline{\Omega}}[r(x)-v_e(x)]$.
  Then there exists unique strictly positive function $N^*\in C(\overline{\Omega})$ such that
\begin{equation}
  HN^*=s(H)N^* \quad \mbox{and} \quad \int_{\Omega}N^*(x)\rmd\mu(x)=1,
\end{equation}
  and any solution of the equation $H\varphi=s(H)\varphi$ is proportional to $N^*$. Moreover,
  $s(H)$ is a dominant spectral value of $H$, i.e., $\mbox{\rm Re}\lambda<s(H)$ for any $\lambda
  \in\sigma(H)\backslash\{s(H)\}$.}
\medskip

  {\em Proof}:\ \ All assertions of this Theorem follow from Lemma 2.6. $\quad\Box$
\medskip

  In the next theorem we consider long-term behavior of the solution $N=N(x,t)$ of the equation
  (2.25) subject to the initial condition $N(\cdot,0)=N_0$, where $N_0\in C(\overline{\Omega})$ is
  a given nonnegative function. Let $H$ and $s(H)$ be as before, and $\rmd\mu\in\mathcal{M}(\Omega)$
  be the positive probability measure in $\Omega$ ensured by the assertion (2) of Theorem 2.7. Set
$$
  M(t)=\displaystyle\int_{\Omega}N(x,t)\rmd\mu(x)\;\; \mbox{for}\;\;t\geq 0,
  \quad \mbox{and} \quad
  M_0=\displaystyle\int_{\Omega}N_0(x)\rmd\mu(x).
$$

  {\bf Theorem 2.8}\ \ {\em  Let notations be as above. In the general case we have
\begin{equation}
  M(t)=M_0e^{s(H)\,t} \quad \mbox{for}\;\; t\geq 0,
\end{equation}
  so that $\displaystyle\lim_{t\to\infty}M(t)=\infty$ if $s(H)>0$, $M(t)=M_0$ for all $t\geq 0$ if
  $s(H)=0$, and $\displaystyle\lim_{t\to\infty}M(t)=0$ if $s(H)<0$. If further the population
  system is ergodic and $s(H)>\displaystyle\max_{x\in\overline{\Omega}}[r(x)-v_e(x)]$ then by
  letting $N^*\in C(\overline{\Omega})$ be the strictly positve function ensured by the assertion
  $(3)$ of Theorem 2.7, we further have
\begin{equation}
  N(x,t)=M_0N^*(x)e^{s(H)\, t}+h(x,t)e^{s(H)\, t} \quad \mbox{for}\;\; x\in\Omega,\;\, t\geq 0,
\end{equation}
  where $h$ is a function uniquely determined by $N_0$ possessing the following properties:

  $(1)$\ \ $\displaystyle\int_{\Omega}h(x,t)\rmd\mu(x)=0$ for $t\geq 0$;

  $(2)$\ \ $\displaystyle\sup_{x\in\Omega}|h(x,t)|\leq Ce^{-\kappa t}\sup_{x\in\Omega}
  |N_0(x)-M_0N^*(x)|$ for $t\geq 0$,

\noindent
  where $\kappa$ and $C$ are positive constants depending only on the migration function $v$.}
\medskip

  {\em Proof}:\ \ First we note that in the general case we have
$$
\begin{array}{rl}
  \displaystyle\frac{\rmd M(t)}{\rmd t}
  =&\displaystyle\int_{\Omega}\frac{\partial N(x,t)}{\partial t}\rmd\mu(x)
  =\int_{\Omega}HN(x,t)\rmd\mu(x)=\langle N(\cdot,t),H'\rmd\mu \rangle
  =\langle N(\cdot,t),s(H)\rmd\mu \rangle
\\[0.3cm]
  =&s(H)\langle N(\cdot,t),\rmd\mu \rangle=s(H)M(t) \quad \mbox{for}\;\; t>0.
\end{array}
$$
  Hence (2.32) holds. Next we assume that the population system is ergodic and $s(H)>\displaystyle
  \max_{x\in\overline{\Omega}}[r(x)-v_e(x)]$. Without loss of generality we may assume that $s(H)=
  0$, for otherwise by replacing $N(x,t)$ and $H$ with $\tilde{N}(x,t)=N(x,t)e^{-s(H)\,t}$ and
  $\tilde{H}=H-s(H)I$, respectively, we reduce the problem into this simpler situation. We denote
$$
  X=C(\overline{\Omega}), \quad X_0=\mbox{\rm span}\{N^*\}, \quad
  X_1=\Big\{\varphi\in X:\int_{\Omega}\varphi(x)\rmd\mu(x)=0\Big\}.
$$
  It is clear that both $X_0$ and $X_1$ are closed linear subspaces of $X$, and $X=X_0\oplus X_1$.
  The assertion (3) of Theorem 2.7 ensures that $X_0=\ker(H)$. Moreover, a similar argument as
  in the proof of (2.32) shows that $\mbox{rg}(H)\subseteq X_1$. Hence, both $X_0$ and $X_1$ are
  invariant subspaces of $H$. It follows that
\begin{equation}
  \sigma(H)=\sigma(H|_{X_0})\cup\sigma(H|_{X_1})=\{0\}\cup\sigma(H|_{X_1}),
\end{equation}
  and
\begin{equation*}
  e^{tH}=e^{tH|_{X_0}}\oplus e^{tH|_{X_1}}=id_{X_0}\oplus e^{tH|_{X_1}} \quad \mbox{for}\;\;  t\geq 0.
\end{equation*}
  We write $N_0(x)=M_0N^*(x)+h_0(x)$, $x\in\Omega$, where $h_0(x)=N_0(x)-M_0N^*(x)$, $x\in\Omega$.
  Then $M_0N^*\in X_0$ and $h_0\in X_1$. It follows by (2.34) that
$$
  N(\cdot,t)=e^{tH}N_0=M_0N^*(x)+e^{tH|_{X_1}}h_0 \quad \mbox{for}\;\;  t\geq 0.
$$
  Letting $h(\cdot,t)=e^{tH|_{X_1}}h_0$, we obtain (2.33). Moreover, since $h(\cdot,t)\in X_1$ for
  any $t\geq 0$, we see that (1) holds. The assertions (2) follows from the fact that $s(H)=0$ is
  a dominant simple eigenvalue of $H$ so that, by (2.34), $s(H|_{X_1})<0$. The proof of Theorem 2.8
  is complete. $\quad\Box$
\medskip

  Combining the assertion (1) of Theorem 2.7 and the first part of Theorem 2.8, we immediately
  obtain the following result:

  {\bf Corollary 2.9}\ \ {\em For a general population system we have the following assertions:

  $(1)$\ \ If the proliferation rate $r$ is everywhere negative, i.e., $r(x)<0$ for all $x\in\Omega$,
  then $\displaystyle\lim_{t\to\infty}M(t)=0$, i.e., the population in this system will finally
  extinct.

  $(2)$\ \ If the proliferation rate $r$ is not everywhere less than or equal to the emigration rate
  $v_e$, i.e., if there exists $x_0\in\Omega$ such that $r(x_0)>v_e(x_0)$, then $\displaystyle
  \lim_{t\to\infty}M(t)=\infty$, i.e., the population in this system will finally become explosive.
  $\quad\Box$}
\medskip

  Another immediate corollary of Theorem 2.8 is that if the population system is ergodic and
  $s(H)>\displaystyle\max_{x\in\overline{\Omega}}[r(x)-v_e(x)]$, then in the case $s(H)=0$ the
  population of this system will finally distributed to the whole habitat domain in a fixed
  proportion $N^*$ uniquely determined by the migration function $v$ of the system, and in the case
  $s(H)>0$ the growth of population in this system exhibits the so-called {\em balanced exponential
  growth} phenomenon, i.e.,
\begin{equation}
  N(x,t)=M_0N^*(x)e^{s(H)\, t}[1+o(1)] \quad \mbox{for}\;\; x\in\Omega, \;\, t\geq 0.
\end{equation}

\subsection{A short discussion to non-ergodic cases}

\hskip 2em
  If we remove the ergodicity assumption then the situation is very much complex. In this
  subsection we only consider two special cases. Moreover, we assume that the population
  system is in the proliferation-stationary case.

  (1) First we consider the case that the habitat domain $\Omega$ is divided into several
  disjoint parts and population in each part forms an independent population system. Hence we
  assume that
\begin{equation*}
  \Omega=\bigcup_{j=1}^m\Omega_j, \quad \mbox{where}\;\;
  \Omega_j\bigcap\Omega_k=\emptyset\;\; \mbox{for}\;\; j\neq k,
\end{equation*}
  and
\begin{equation*}
  v(x,y)=0 \quad \mbox{if}\;\; (x,y)\in \Omega_j\times\Omega_k,\;\; j\neq k,
\end{equation*}
  where $j,k=1,2,\cdots,m$. Moreover, for every $1\leq j\leq m$ we assume that
  $v|_{\Omega_j\times\Omega_j}$ is continuous and ergodic in $\Omega_j\times\Omega_j$. Thus,
  the whole population system is divided into several subsystems, with each subsystem being
  ergodic and different subsystems being mutually independent or having no interchange of
  population between different subsystems.

  By first applying Theorems 2.3 and 2.4 to each subsystem and next combing the results together,
  we easily obtain the following result:
\medskip

  {\bf Theorem 2.10}\ \ {\em Let assumptions be as above. Then the equation $(2.6)$ has a
  solution $N^*$ which is strictly positive and piecewise continuous in $\overline{\Omega}$,
  and is unique under the normalization conditions $\displaystyle\int_{\Omega_j}\!\!N^*(x)
  \rmd x=1$, $j=1,2,\cdots,m$. Moreover, for any $N_0\in C(\overline{\Omega})$, by letting
  $M_{0j}=\displaystyle\int_{\Omega_j}N_0(x)\rmd x$ $(j=1,2,\cdots,m)$ and $N=N(x,t)$ the solution of
  the problem $(2.3)$--$(2.4)$, we have
\begin{equation*}
  N(x,t)=\sum_{j=1}^m [M_{0j}N_j^*(x)+h_j(x,t)] \quad
  \mbox{for}\;\; x\in\Omega,\;\, t\geq 0,
\end{equation*}
  where $N_j^*(x)=N^*(x)\chi_{\Omega_j}(x)$ with $\chi_{\Omega_j}$ being the characteristic function
  of $\Omega_j$ $(j=1,2,\cdots,m)$, and $h_j$'s are functions uniquely determined by
  $N_0$ which possesses the following properties for all  $1\leq j\leq m$:

  $(1)$\ \ $h_j(x,t)=0$ for all $x\in\Omega\backslash\Omega_j$ and $t\geq 0$, and $\displaystyle
  \int_{\Omega_j}h_j(x,t)\rmd x=0$ for all $t\geq 0$;

  $(2)$\ \ $\displaystyle\sup_{x\in\Omega_j}|h_j(x,t)|\leq Ce^{-\kappa t}\sup_{x\in\Omega_j}
  \Big|N_0(x)-M_{0j}N_j^*(x)\Big|$ for all $t\geq 0$;

  $(3)$\ \ $\displaystyle\int_{\Omega_j}|h_j(x,t)|\rmd x\leq Ce^{-\kappa t}\int_{\Omega_j}
  \Big|N_0(x)-M_{0j}N_j^*(x)\Big|\rmd x$ for all $t\geq 0$.

\noindent
  Here $\kappa$ and $C$ are positive constants depending only on the function $v$. $\quad\Box$}
\medskip

  Note that Theorem 2.5 implies that if $N_0\in C(\overline{\Omega})$ then
\begin{equation*}
  \lim_{t\to\infty}\sup_{x\in\Omega}\Big|N(x,t)-\sum_{j=1}^m M_{0j}N_j^*(x)\Big|=0,
\end{equation*}
  and if $N_0\in L^1(\Omega)$ then
\begin{equation*}
  \lim_{t\to\infty}\int_{\Omega}\Big|N(x,t)-\sum_{j=1}^m M_{0j}N_j^*(x)\Big|\rmd x=0.
\end{equation*}
  Moreover, the convergence is exponentially fast.
\medskip

  Comparing the above result with the corresponding one obtained in the Subsection 2.4, we see
  that unlike the ergodic population system, for the present separable population system the
  development of population in the system is not integrated: If $N_0|_{\Omega_j}=0$ for some
  $1\leq j\leq m$ then population in this $\Omega_j$ will be always zero, and development of
  populations in the other habitat domains does not have any influence to this domain.

  (2) Next we consider the case that the habitat domain $\Omega$ is divided into two disjoint
  parts $\Omega_1$ and $\Omega_2$, i.e., $\Omega_1\cap\Omega_2=\emptyset$ and $\Omega=\Omega_1
  \cup\Omega_2$, such that $v(x,y)>0$ for $(x,y)\in\Omega_1\times\Omega_1$, $\Omega_1\times
  \Omega_2$ and $\Omega_2\times\Omega_2$, whereas $v(x,y)=0$ for $(x,y)\in\Omega_2\times\Omega_1$.

  Let $N_j=N|_{\Omega_j\times[0,\infty)}$, $N_{j0}=N_0|_{\Omega_j}$ (for given nonnegative $N_0\in
  C(\overline{\Omega})$), $j=1,2$, and
$$
  v_e^1(x)=\int_{\Omega_1}v(y,x)\rmd y \;\;\; \mbox{for}\;\; x\in\Omega_1, \quad
  v_e^2(x)=\int_{\Omega}v(y,x)\rmd y \;\;\; \mbox{for}\;\; x\in\Omega_2.
$$
  Then the equation (2.3) can be rewritten as the following system of equations:
\begin{equation}
  \frac{\partial N_1(x,t)}{\partial t}=\!\!\int_{\Omega_1}\!\!\!v(x,y)N_1(y,t)\rmd y
  +\!\!\int_{\Omega_2}\!\!\!v(x,y)N_2(y,t)\rmd y-v_e^1(x)N_1(x,t), \quad
  x\in\Omega_1,\;\; t>0,
\end{equation}
\begin{equation}
  \frac{\partial N_2(x,t)}{\partial t}=\int_{\Omega_2}\!\!v(x,y)N_2(y,t)\rmd y-v_e^2(x)N_2(x,t),
  \quad  x\in\Omega_2,\;\; t>0.
\end{equation}
  Clearly, the equation (2.37) can be decoupled from the equation (2.36). Let
$$
  v_e^{21}(x)=\int_{\Omega_1}v(y,x)\rmd y \quad \mbox{and}  \quad
  v_e^{22}(x)=\int_{\Omega_2}v(y,x)\rmd y \;\;\; \mbox{for}\;\; x\in\Omega_2.
$$
  Then $v_e^2=v_e^{21}+v_e^{22}$ and $v_e^{21}$, $v_e^{22}$ are strictly positive in $\Omega_2$.
  We rewrite the equation (2.37) as follows:
\begin{equation*}
  \frac{\partial N_2(x,t)}{\partial t}=\int_{\Omega_2}\!\!v(x,y)N_2(y,t)\rmd y-v_e^{22}(x)N_2(x,t)
  -v_e^{21}(x)N_2(x,t), \quad  x\in\Omega_2,\;\; t>0.
\end{equation*}
  This equation can be regarded as a population migration equation in the domain $\Omega_2$ with
  proliferation rate $r(x)=-v_e^{21}(x)$. Since the proliferation rate is strictly negative
  everywhere in $\Omega_2$, by the assertion (1) of Corollary 2.9 we see that $N_2(x,t)$ converges
  to zero in $\Omega_2$ as $t\to\infty$. More precisely, by using the assertion (1) of
  Theorem 2.7 and the expression (2.33) in Theorem 2.8 we see that there exists a positive
  constant $\lambda_2\geq\displaystyle\inf_{x\in\Omega_2}v_e^{21}(x)$ such that
\begin{equation}
  N_2(x,t)=M_{20}N_2^*(x)e^{-\lambda_2 t}+h_2(x,t)e^{-\lambda_2\, t} \quad
  \mbox{for}\;\; x\in\Omega_2,\;\, t\geq 0,
\end{equation}
  where $M_{20}=\displaystyle\!\int_{\Omega_2}\!\!N_0(x)\rmd\mu_2(x)$ with $\rmd\mu_2$ being a
  positive probability measure in $\Omega_2$, $N_2^*$ is a strictly positive continuous function
  in $\overline{\Omega}_2$ such that $\displaystyle\!\int_{\Omega_2}\!\!N_2^*(x)\rmd\mu_2(x)=1$,
  and $h_2$ is a function defined in $\overline{\Omega}_2\times [0,\infty)$ satisfying the
  following condition:
\begin{equation*}
  \sup_{\Omega_2}|h_2(x,t)|\leq Ce^{-\kappa_2 t}\sup_{\Omega_2}|N_2(x,0)-M_{20}N_2^*(x)| \quad
  \mbox{for}\;\; x\in\Omega_2,\;\, t\geq 0,
\end{equation*}
  where $\kappa_2$ and $C$ are positive constants depending only on $v$. By (2.38) we see that
  as $t\to\infty$, $N_2(x,t)$ uniformly converges to zero in $\Omega_2$ with an exponential
  speed.

  On the other hand, since the total amount of the population in the whole domain $\Omega$, i.e.,
  the quantity $M(t)=\displaystyle\!\int_{\Omega}\!N(x,t)\rmd x$, is constant or identically equal
  to $M_0=\displaystyle\!\int_{\Omega}\!N_0(x)\rmd x$ for all $t\geq 0$, we see that the whole
  population that might be originally distributed in the whole domain $\Omega$ will finally
  distributed to the subdomain $\Omega_1$. The explicit expression of $N_1(x,t)$ can be easily
  obtained by substituting the expression (2.38) into (2.36) and next solving that equation by using
  Theorem 2.5 to the semigroup $e^{tH_1}$ in $C(\overline{\Omega}_1)$, where $H_1$ is the bounded
  linear operator in $C(\overline{\Omega}_1)$ defined by
$$
  H_1\varphi(x)=\int_{\Omega_1}\!\!\!v(x,y)\varphi(y)\rmd y-v_e^1(x)\varphi(x) \quad
   \mbox{for}\;\; \varphi\in C(\overline{\Omega}_1).
$$
  We omit the details. Note that the equation (2.36) can be regarded as a proliferation-stationary
  population evolution equation with a source term; the second term on the right-hand side of that
  equation is the source term.

\section{Migration epidemics models}
\setcounter{equation}{0}

\hskip 2em
  In this section we consider three simple mathematical models describing spread of epidemics
  in migration population systems. These models include: The migration SI model, the migration
  SIR model, and the migration SIRE model. Limited to spaces, we shall be content with deducing
  some elementary results, and leave more penetrating study for future work.

\subsection{The migration SI model}

\hskip 2em
  We first consider the migration SI model, which is obtained by adding migration terms into the
  classical SI model. Let $S(x,t)$ be the density of the susceptible population, and
  $I(x,t)$ be the density of the infective population. Then we have the following system
  of equations:
\begin{equation}
\left\{
\begin{array}{l}
  \displaystyle\frac{\partial S(x,t)}{\partial t}=\int_{\Omega}v(x,y)S(y,t)\rmd y-v_e(x)S(x,t)
  -rS(x,t)I(x,t), \quad x\in\Omega,\;\; t>0,
\\ [0.3cm]
  \displaystyle\frac{\partial I(x,t)}{\partial t}=\int_{\Omega}v(x,y)I(y,t)\rmd y-v_e(x)I(x,t)
  +rS(x,t)I(x,t), \quad x\in\Omega,\;\; t>0.
\end{array}
\right.
\end{equation}
  We impose this system with the following initial condition:
\begin{equation}
  S(x,0)=S_0(x), \quad I(x,0)=I_0(x) \quad \mbox{for}\;\; x\in\Omega,
\end{equation}
  where $S_0$ and $I_0$ are given functions.

  Unlike (2.4) which is a linear equation, the system of equations (2.14) are nonlinear equations.
  Hence, unlike the initial value problem (2.3)--(2.4), for the initial value problem (3.1)--(3.2)
  both the standard Picard iteration argument and the $C_0$-semigroup theory only ensure local
  existence of a solution. More precisely, for any $S_0,I_0\in C(\overline{\Omega})$, there exists
  corresponding $\delta>0$ such that the initial value problem (2.3)--(2.4) has a unique solution
  $S,I\in C(\overline{\Omega}\times[0,\delta])$. By a standard argument, the solution can be uniquely
  extended into a maximal time interval $[0,T)$, where $0<T\leq\infty$, and in case $T<\infty$ the
  solution $(S,I)$ must blow-up as $t\to T^-$. As we shall see below, for nonnegative initial
  values $S_0$ and $I_0$, the solution of the above problem satisfies
\begin{equation}
  0\leq S(x,t)\leq C, \quad  0\leq I(x,t)\leq C
\end{equation}
  for all $(x,t)$ in the existence domain of the solution, where $C$ is a constant depending on
  $(S_0,I_0)$. It follows that for nonnegative initial values the solution of the above problem
  exists globally.

  The assertion that $S(x,t)\geq 0$ and $I(x,t)\geq 0$ for all $(x,t)$ in the existence region
  of the solution follows from Lemma 0.0 in the appendix. Let
\begin{equation}
  N(x,t)=S(x,t)+I(x,t) \quad \mbox{and} \quad  N_0(x)=S_0(x)+I_0(x).
\end{equation}
  From (3.1) and (3.2) it can be easily seen that $N=N(x,t)$ satisfies (2.4) and (2.5). It follows
  that $N=N(x,t)$ has the following expression:
\begin{equation}
  N(x,t)=M_0N^*(x)+h(x,t) \quad \mbox{for}\;\; x\in\Omega,\;\, t\geq 0,
\end{equation}
  where $h$ satisfies the following properties:
\begin{equation}
  \int_{\Omega}h(x,t)\rmd x=0, \quad
  \sup_{x\in\Omega}|h(x,t)|\leq C\sup_{x\in\Omega}|N_0(x)-M_0N^*(x)|e^{-\kappa t}
  \quad \mbox{for}\;\; t\geq 0,
\end{equation}
  where $\kappa$ and $C$ are positive constants determined by the migration function $v$ (so that
  it is independent of the initial data). This particularly implies that $N$ is a bounded function.
  Since (3.4) implies that $0\leq S(x,t)\leq N(x,t)$ and $0\leq I(x,t)\leq N(x,t)$, we see that
  (3.3) follows.

  From (3.1) we easily see that
$$
  \frac{\rmd }{\rmd t}\Big(\int_{\Omega}S(x,t)\rmd x\Big)\leq 0 \quad \mbox{and} \quad
  \frac{\rmd }{\rmd t}\Big(\int_{\Omega}I(x,t)\rmd x\Big)\geq 0.
$$
  This means that the total amount of susceptible population is decreasing, whereas the total amount
  of infective population is increasing.

  Before considering the asymptotic behavior of solutions of the problem (3.1) and (3.2), let us
  first consider the steady-state solution of the equation (3.1). Let $(S(x),I(x))$ be a
  steady-state solution of the equation (3.1). Clearly $(S(x),I(x))$ satisfies the following system
  of equations:
\begin{equation}
\left\{
\begin{array}{l}
  \displaystyle\int_{\Omega}v(x,y)S(y)\rmd y-v_e(x)S(x)-rS(x)I(x)=0,
   \quad x\in\Omega,
\\ [0.3cm]
  \displaystyle\int_{\Omega}v(x,y)I(y)\rmd y-v_e(x)I(x)+rS(x)I(x)=0,
   \quad x\in\Omega.
\end{array}
\right.
\end{equation}
  Let $N(x)=S(x)+I(x)$. Clearly $N=N(x)$ satisfies the equation (2.6). Hence there exists a
  constant $c$ such that $N(x)=cN^*(x)$, or in other words,
$$
  S(x)+I(x)=cN^*(x) \quad \mbox{for}\;\; x\in\Omega.
$$
  Integration either one of the two equations in (3.5), we get
$$
  \int_{\Omega}S(x)I(x)\rmd x=0.
$$
  Since we are only considering nonnegative solutions, this implies that
$$
  S(x)I(x)=0 \quad \mbox{for}\;\, x\in\Omega.
$$
  so that one of the following three situations occurs: $(i)$ $(S(x),I(x))=(cN^*(x),0)$;
  $(ii)$ $(S(x),I(x))=(0,cN^*(x))$; $(iii)$ there exist disjoint measurable subsets $\Omega_1$
  and $\Omega_2$ of $\Omega$ with positive measures, such that $\Omega=\Omega_1\cup\Omega_2$
  and $(S(x),I(x))=(cN^*(x)\chi_{\Omega_1}(x),cN^*(x)\chi_{\Omega_2}(x))$. The third situation
  cannot occur. Indeed, otherwise from (3.5) we get
\begin{equation*}
\left\{
\begin{array}{l}
  \displaystyle\int_{\Omega_1}v(x,y)N^*(y)\rmd y=v_e(x)N^*(x)\chi_{\Omega_1}(x),
   \quad x\in\Omega,
\\ [0.3cm]
  \displaystyle\int_{\Omega_2}v(x,y)N^*(y)\rmd y=v_e(x)N^*(x)\chi_{\Omega_2}(x),
   \quad x\in\Omega,
\end{array}
\right.
\end{equation*}
  which implies that
\begin{equation*}
  \int_{\Omega_1}v_{o2}(y)N^*(y)\rmd y=0, \qquad
  \int_{\Omega_2}v_{o1}(y)N^*(y)\rmd y=0,
\end{equation*}
  where $v_{o2}(y)=\displaystyle\int_{\Omega_2}v(x,y)\rmd x$, $v_{o1}(y)=\displaystyle\int_{\Omega_1}
  v(x,y)\rmd x$. Since $v_{o2}(y)>0$, $v_{o1}(y)>0$, we conclude that $N^*(y)=0$ for all $x\in\Omega$,
  which is a contradiction. This proves the desired assertion. Hence, (3.7) has only two
  linearly independent nontrivial nonnegative solutions; they are respectively $(S(x),I(x))=
  (cN^*(x),0)$ and $(S(x),I(x))=(0,cN^*(x))$, where $c$ is a positive constant.

  As for asymptotic behavior of solutions of the problem (3.1) and (3.2), we have the following
  result:
\medskip

  {\bf Theorem 3.1}\ \ {\em Assume that $v\in C(\overline{\Omega}\times\overline{\Omega})$
  and $v(x,y)>0$ for all $x,y\in\overline{\Omega}$. For $S_0, I_0\in C(\overline{\Omega})$
  such that $S_0\geq0$ and $I_0{\geq\atop\neq}0$, let $N_0=S_0+I_0$ and $M_0=\displaystyle
  \int_{\Omega}N_0(x)\rmd x$. For the solution $(S,I)$ of the problem $(3.1)$--$(3.2)$ we have the
  following relation:
\begin{equation}
  \lim_{t\to\infty}(S(x,t),I(x,t))=(0,M_0N^*(x)) \quad \mbox{uniformly for}\;\; x\in\Omega.
\end{equation}
}

  {\em Proof}:\ \ Let $\underline{I}=\underline{I}(x,t)$ be the solution of the following
  initial value problem:
\begin{equation*}
\left\{
\begin{array}{l}
  \displaystyle\frac{\partial\underline{I}(x,t)}{\partial t}
  =\int_{\Omega}v(x,y)\underline{I}(y,t)\rmd y-v_e(x)\underline{I}(x,t), \quad x\in\Omega,\;\; t>0,
\\ [0.1cm]
  \underline{I}(x,0)=I_0(x), \quad x\in\Omega.
\end{array}
\right.
\end{equation*}
  It is easy to deduce that
\begin{equation*}
  I(x,t)\geq\underline{I}(x,t) \quad \quad \mbox{for}\;\; x\in\Omega,\;\, t\geq 0.
\end{equation*}
  By Theorem 2.5 we have
\begin{equation*}
  \lim_{t\to\infty}\underline{I}(x,t)=cN^*(x) \quad \mbox{uniformly for}\;\; x\in\Omega,
\end{equation*}
  where $c=\displaystyle\int_{\Omega}I_0(x)\rmd x>0$. Since $\displaystyle\inf_{x\in\Omega}
  N^*(x)>0$, it follows that there exist positive constants $T$ and $\delta$ such that
\begin{equation*}
  I(x,t)\geq\delta \quad \quad \mbox{for}\;\; x\in\Omega,\;\, t\geq T.
\end{equation*}
  Hence, from the equation satisfied by $S$ we get
\begin{equation*}
  \frac{\partial S(x,t)}{\partial t}\leq\int_{\Omega}v(x,y)S(y,t)\rmd y-v_e(x)S(x,t)
  -r\delta S(x,t), \quad x\in\Omega,\;\; t\geq T,
\end{equation*}
  which is equivalent to
\begin{equation*}
  \frac{\partial}{\partial t}\Big(S(x,t)e^{r\delta t}\Big)\leq
  \int_{\Omega}v(x,y)S(y,t)e^{r\delta t}\rmd y-v_e(x)S(x,t)e^{r\delta t},
   \quad x\in\Omega,\;\; t\geq T.
\end{equation*}
  Let $u=u(x,t)$ ($x\in\Omega$, $t\geq T$) be the solution of the following initial value problem:
\begin{equation*}
\left\{
\begin{array}{l}
  \displaystyle\frac{\partial u(x,t)}{\partial t}
  =\int_{\Omega}v(x,y)u(y,t)\rmd y-v_e(x)u(x,t), \quad x\in\Omega,\;\; t>T,
\\ [0.1cm]
  u(x,T)=S(x,T)e^{r\delta T}, \quad x\in\Omega.
\end{array}
\right.
\end{equation*}
  It follows that
\begin{equation*}
  S(x,t)e^{r\delta t}\leq u(x,t) \quad \quad \mbox{for}\;\; x\in\Omega,\;\, t\geq T.
\end{equation*}
  Since by Theorem 2.5 we know that $u=u(x,t)$ ($x\in\Omega$, $t\geq T$) is a bounded function,
  the above inequality implies that
\begin{equation*}
  \lim_{t\to\infty}S(x,t)=0 \quad \mbox{uniformly for}\;\; x\in\Omega,
\end{equation*}
  and, consequently,
\begin{equation*}
  \lim_{t\to\infty}I(x,t)=\lim_{t\to\infty}[N(x,t)-S(x,t)]=M_0N^*(x)
  \quad \mbox{uniformly for}\;\; x\in\Omega.
\end{equation*}
  This proves Theorem 3.1. $\quad\Box$
\medskip

\subsection{The migration SIR model}

\hskip 2em
  We Next consider the migration SIR model, which is obtained by adding migration terms into the
  classical SIR model. Let $S(x,t)$ be the density of the susceptible population, $I(x,t)$ be the
  density of the infective population, and $R(x,t)$ be the density of the recovered population.
  Then we have the following system of equations:
\begin{equation}
\left\{
\begin{array}{l}
  \displaystyle\frac{\partial S(x,t)}{\partial t}=\int_{\Omega}v(x,y)S(y,t)\rmd y-v_e(x)S(x,t)
  -rS(x,t)I(x,t), \quad x\in\Omega,\;\; t>0,
\\ [0.3cm]
  \displaystyle\frac{\partial I(x,t)}{\partial t}=\int_{\Omega}v(x,y)I(y,t)\rmd y-v_e(x)I(x,t)
  +rS(x,t)I(x,t)-aI(x,t), \quad x\in\Omega,\;\; t>0,
\\ [0.3cm]
  \displaystyle\frac{\partial R(x,t)}{\partial t}=\int_{\Omega}v(x,y)R(y,t)\rmd y-v_e(x)R(x,t)
  +aI(x,t), \quad x\in\Omega,\;\; t>0.
\end{array}
\right.
\end{equation}
  We impose this system with the following initial condition:
\begin{equation}
  S(x,0)=S_0(x), \quad I(x,0)=I_0(x), \quad R(x,0)=R_0(x) \quad \mbox{for}\;\; x\in\Omega,
\end{equation}
  where $S_0$, $I_0$ and $R_0$ are given functions.

  Note that (3.9) is also a system of nonlinear differential-equations. By a similar argument as
  in the subsection 3.1, we can prove that for any $S_0,I_0,R_0\in C(\overline{\Omega})$ such that
  $S_0,I_0,R_0\geq 0$, the initial value problem (3.9)--(3.10) has a unique global solution
  $S,I,R\in C(\overline{\Omega}\times[0,\infty))$, such that
\begin{equation}
  S(x,t)\geq 0, \quad  I(x,t)\geq 0, \quad  R(x,t)\geq 0
\end{equation}
  for all $x\in\Omega$ and $t\geq 0$. Moreover, by letting
$$
  N(x,t)=S(x,t)+I(x,t)+R(x,t) \quad \mbox{and} \quad N_0(x)=S_0(x)+I_0(x)+R_0(x),
$$
  we have that $N=N(x,t)$ is a solution of the problem (2.4) and (2.5), so that $N=N(x,t)$
  satisfies (3.5) and (3.6).

  From the first and the last equations in (3.9) we easily see that
$$
  \frac{\rmd }{\rmd t}\Big(\int_{\Omega}S(x,t)\rmd x\Big)\leq 0 \quad \mbox{and} \quad
  \frac{\rmd }{\rmd t}\Big(\int_{\Omega}R(x,t)\rmd x\Big)\geq 0.
$$
  This means that the amount of susceptible population is decreasing, and the amount of recovered
  population is increasing.

  Let us now study asymptotic behavior of solutions of the problem (3.9) and (3.10). We first
  consider the steady-state solution of the equation (3.9). Let $(S(x),I(x),R(x))$ be a
  steady-state solution of the equation (3.9). Clearly $(S(x),I(x),R(x))$ satisfies the following
  system of equations:
\begin{equation}
\left\{
\begin{array}{l}
  \displaystyle\int_{\Omega}v(x,y)S(y)\rmd y-v_e(x)S(x)-rS(x)I(x)=0,
   \quad x\in\Omega,
\\ [0.3cm]
  \displaystyle\int_{\Omega}v(x,y)I(y)\rmd y-v_e(x)I(x)+rS(x)I(x)-aI(x)=0,
   \quad x\in\Omega,
\\ [0.3cm]
  \displaystyle\int_{\Omega}v(x,y)R(y)\rmd y-v_e(x)R(x)+aI(x)=0,
   \quad x\in\Omega.
\end{array}
\right.
\end{equation}
  Let $N(x)=S(x)+I(x)+R(x)$. Clearly $N=N(x)$ satisfies the equation (2.6). Hence there exists a
  constant $c$ such that $N(x)=cN^*(x)$, or in other words,
$$
  S(x)+I(x)+R(x)=cN^*(x) \quad \mbox{for}\;\; x\in\Omega.
$$
  Integration the last equation in (3.5), we get
$$
  \int_{\Omega}I(x)\rmd x=0.
$$
  Hence $I(x)=0$ for $x\in\Omega$. Substituting this relation into the first and the last equations
  in (3.12), we get
$$
  \int_{\Omega}v(x,y)S(y)\rmd y-v_e(x)S(x)=0 \quad \mbox{and} \quad
  \int_{\Omega}v(x,y)R(y)\rmd y-v_e(x)R(x)=0.
$$
  By Theorem 2.1, this implies that $S(x)=c_1N^*(x)$ and $R(x)=c_2N^*(x)$ for $x\in\Omega$.
  Hence, steady-state solution of the equation (3.9) has the form $(S(x),I(x),R(x))=
  (c_1N^*(x),0,c_2N^*(x))$, where $c_1$ and $c_2$ are nonnegative constants.

  We denote by $M_S(t)$, $M_I(t)$ and $M_R(t)$ the amounts of susceptible, infective, and
  recovered populations, respectively, at time $t$. In the above we have seen that $M_S(t)$
  and $M_R(t)$ are respectively monotone decreasing and monotone increasing, so that they both
  have limits as $t\to\infty$. Since $M_S(t)+M_I(t)+M_R(t)=M(t)=M_0$, where $M(t)$ and $M_0$
  denote the total amount of the population at time $t$ and the initial total amount of the
  population, respectively, we infer that $M_I(t)$ also has limit as $t\to\infty$. We now prove
  the following result:
\medskip

  {\bf Theorem 3.2}\ \ {\em Assume that $v\in C(\overline{\Omega}\times\overline{\Omega})$
  and $v(x,y)>0$ for all $x,y\in\overline{\Omega}$. For $S_0,I_0,R_0\in C(\overline{\Omega})$
  such that $I_0\geq0$, $R_0\geq0$ and $S_0{\geq\atop\neq}0$, let $N_0=S_0+I_0+R_0$ and $M_0=
  \displaystyle\int_{\Omega}N_0(x)\rmd x$. For the solution $(S,I,R)$ of the problem $(3.9)$--$(3.10)$
  we have the following relations:
\begin{equation}
\left\{
\begin{array}{l}
  S(x,t)=M_S(t)N^*(x)+h_S(x,t),\\
  I(x,t)=M_I(t)N^*(x)+h_I(x,t),\\
  R(x,t)=M_R(t)N^*(x)+h_R(x,t),
\end{array}
\right.
\end{equation}
  where $\displaystyle\lim_{t\to\infty}(M_S(t),M_I(t),M_R(t))=(c_1,0,c_2)$ for some positive
  constants $c_1$, $c_2$ such that $c_1+c_2=M_0$, and for $i=S,I,R$,
\begin{equation}
  \int_{\Omega}h_i(x,t)\rmd x=0, \;\; \forall t\geq 0 \quad \mbox{and} \quad
  \lim_{t\to\infty}\int_{\Omega}|h_i(x,t)|\rmd x=0
\end{equation}
}

  {\em Proof}:\ \ From the third equation in (3.9) we have
$$
  \frac{\rmd }{\rmd t}\Big(\int_{\Omega}R(x,t)\rmd x\Big)=a\int_{\Omega}I(x,t)\rmd x
  \quad \mbox{for all}\;\, t>0.
$$
  Since both $M_S(t)=\displaystyle\int_{\Omega}S(x,t)\rmd x$ and $M_R(t)=\displaystyle\int_{\Omega}
  R(x,t)\rmd x$ have limits as $t\to\infty$, from the above relation we infer that there must hold
  the relation
$$
  \lim_{t\to\infty}M_I(t)=\lim_{t\to\infty}\int_{\Omega}I(x,t)\rmd x=0,
$$
  for otherwise a contradiction will follow. Since
$$
  \int_{\Omega}|h_I(x,t)|\rmd x=\int_{\Omega}|I(x,t)-M_I(t)N^*(x)|\rmd x
  \leq\int_{\Omega}I(x,t)\rmd x+M_I(t)=2\int_{\Omega}I(x,t)\rmd x,
$$
  we conclude that
$$
  \lim_{t\to\infty}\int_{\Omega}|h_I(x,t)|\rmd x=0.
$$
  Next, let
$$
  f(x,t)=-r\Big[S(x,t)I(x,t)-\Big(\int_{\Omega}S(x,t)I(x,t)\rmd x\Big)N^*(x)\Big].
$$
  From the first equation in (3.9) we get
$$
  \frac{\partial h_S(x,t)}{\partial t}=\int_{\Omega}v(x,y)h_S(y,t)\rmd y-v_e(x)h_S(x,t)
  +f(x,t), \quad x\in\Omega,\;\; t>0.
$$
  Since $\displaystyle\int_{\Omega}f(x,t)\rmd x=0$ and the closed linear subspace $X_1=\Big\{\varphi
  \in C(\overline{\Omega}):\displaystyle\int_{\Omega}\varphi(x)\rmd x=0\Big\}$ of $X=
  C(\overline{\Omega})$ is an invariant subspace of the bounded linear operator $H:X\to
  X$ defined by
$$
  H(\varphi)=\int_{\Omega}v(x,y)\varphi(y)\rmd y-v_e(x)\varphi(x)
  \quad \mbox{for} \;\, \varphi\in X,
$$
  the above equation can be regarded as a differential equation in this subspace $X_1$. More
  precisely, denoting by $H_1$ the restriction of $H$ in $X_1$, the above
  equation can be rewritten in the following form:
$$
  \frac{\rmd  h_S}{\rmd t}=H_1h_S+f.
$$
  Hence its solution can be expressed as follows:
$$
  h_S(\cdot,t)=e^{tH_1}h_S(\cdot,0)+\int_0^te^{(t-s)H_1}f(\cdot,s)ds.
$$
  By the assertion (3) of Theorem 2.5 we thus get
$$
  \|h_S(\cdot,t)\|_{L^1(\Omega)}\leq C\|h_S(\cdot,0)\|_{L^1(\Omega)}e^{-\kappa t}
  +C\int_0^t\|f(\cdot,s)\|_{L^1(\Omega)}e^{-\kappa(t-s)}ds.
$$
  Since $S(x,t)$ is a bounded function and $\displaystyle\lim_{t\to\infty}\int_{\Omega}I(x,t)\rmd x=0$,
  we see that
$$
  \lim_{t\to\infty}\|f(\cdot,t)\|_{L^1(\Omega)}=0.
$$
  From the above two relations we easily deduce that
$$
  \lim_{t\to\infty}\int_{\Omega}|h_S(x,t)|\rmd x=\lim_{t\to\infty}\|h_S(\cdot,t)\|_{L^1(\Omega)}=0.
$$
  By a similar argument we can also prove that
$$
  \lim_{t\to\infty}\int_{\Omega}|h_R(x,t)|\rmd x=\lim_{t\to\infty}\|h_R(\cdot,t)\|_{L^1(\Omega)}=0.
$$
  This proves Theorem 3.2. $\quad\Box$
\medskip

\subsection{The migration SIRE model}

\hskip 2em
  Finally let us consider the migration SIRE model, which is obtained by adding migration terms into
  the classical SIRE model. Let $S(x,t)$ be the density of the susceptible population, $I(x,t)$ be
  the density of the infective population, and $R(x,t)$ be the density of the recovered population.
  Then we have the following system of equations:
\begin{equation}
\left\{
\begin{array}{l}
  \displaystyle\frac{\partial S(x,t)}{\partial t}=\int_{\Omega}v(x,y)S(y,t)\rmd y-v_e(x)S(x,t)
  -rS(x,t)I(x,t)+bR(x,t), \quad x\in\Omega,\;\; t>0,
\\ [0.3cm]
  \displaystyle\frac{\partial I(x,t)}{\partial t}=\int_{\Omega}v(x,y)I(y,t)\rmd y-v_e(x)I(x,t)
  +rS(x,t)I(x,t)-aI(x,t), \quad x\in\Omega,\;\; t>0,
\\ [0.3cm]
  \displaystyle\frac{\partial R(x,t)}{\partial t}=\int_{\Omega}v(x,y)R(y,t)\rmd y-v_e(x)R(x,t)
  +aI(x,t)-bR(x,t), \quad x\in\Omega,\;\; t>0.
\end{array}
\right.
\end{equation}
  We impose it with the following initial conditions:
\begin{equation}
  S(x,0)=S_0(x), \quad I(x,0)=I_0(x), \quad R(x,0)=R_0(x) \quad \mbox{for}\;\; x\in\Omega,
\end{equation}
  where $S_0$, $I_0$ and $R_0$ are given functions.

  Note that (3.15) is also a system of nonlinear differential-equations. By a similar argument as
  in the subsection 3.1, we can prove that for any $S_0,I_0,R_0\in C(\overline{\Omega})$ such that
  $S_0,I_0,R_0\geq 0$, the initial value problem (3.15)--(3.16) has a unique global solution
  $S,I,R\in C(\overline{\Omega}\times[0,\infty))$, such that
\begin{equation}
  S(x,t)\geq 0, \quad  I(x,t)\geq 0, \quad  R(x,t)\geq 0
\end{equation}
  for all $x\in\Omega$ and $t\geq 0$. Moreover, by letting
$$
  N(x,t)=S(x,t)+I(x,t)+R(x,t) \quad \mbox{and} \quad N_0(x)=S_0(x)+I_0(x)+R_0(x),
$$
  we have that $N=N(x,t)$ is a solution of the problem (2.4) and (2.5), so that $N=N(x,t)$
  satisfies (3.5) and (3.6). from the above relation we see that $S(x,t)=N(x,t)-I(x,t)-R(x,t)$.
  Substituting this relation into the second equation in (3.15) and coupling it with the third
  equation, we obtain the following reduced system:
\begin{equation}
\left\{
\begin{array}{l}
  \displaystyle\frac{\partial I(x,t)}{\partial t}=\int_{\Omega}v(x,y)I(y,t)\rmd y
  -[v_e(x)+a-rN(x,t)]I(x,t)-rI^2(x,t)-rI(x,t)R(x,t),
\\
   \quad\quad\quad\quad\quad\quad\quad\quad\quad\quad x\in\Omega,\;\; t>0,
\\ [0.3cm]
  \displaystyle\frac{\partial R(x,t)}{\partial t}=\int_{\Omega}v(x,y)R(y,t)\rmd y-v_e(x)R(x,t)
  +aI(x,t)-bR(x,t), \quad x\in\Omega,\;\; t>0.
\end{array}
\right.
\end{equation}

  A steady-state solution of the system (3.18) is a solution of the following system:
\begin{equation}
\left\{
\begin{array}{l}
  \displaystyle\int_{\Omega}v(x,y)I(y)\rmd y-[v_e(x)+a-rN(x)]I(x)-rI^2(x)-rI(x)R(x)=0,
   \quad x\in\Omega,
\\ [0.3cm]
  \displaystyle\int_{\Omega}v(x,y)R(y)\rmd y-[v_e(x)+b]R(x)+aI(x)=0, \quad x\in\Omega,
\end{array}
\right.
\end{equation}
  where $N(x)$ is a steady-state solution of (2.4), or a solution of (2.6), so that $N(x)=M_0
  N^*(x)$. From the second equation in (3.19) we have
$$
  R=aL_bI,
$$
  where $L_b$ is a positive bounded linear operator in $C(\overline{\Omega})$. Substituting
  this expression as well as the relation $N(x)=M_0N^*(x)$ into the first equation in (3.19), we
  get the following equation:
\begin{equation}
  \int_{\Omega}v(x,y)I(y)\rmd y-[v_e(x)+a-rM_0N^*(x)]I(x)-rI(x)(1+aL_b)I(x)=0,
   \quad x\in\Omega.
\end{equation}
   The linearization of the above equation at the trivial solution $I=0$ is as follows:
\begin{equation}
  \int_{\Omega}v(x,y)\varphi(y)\rmd y-[v_e(x)+a-rM_0N^*(x)]\varphi(x)=0,
   \quad x\in\Omega.
\end{equation}

  For every $r\in\mathbb{R}$, let $A_r$ be the bounded linear operator in $C(\overline{\Omega})$
  defined as follows: For any $\varphi\in C(\overline{\Omega})$,
$$
  A_r\varphi=\mbox{the left-hand side of (3.21)}.
$$
  It is clear that the uniformly continuous semigroup generated by $A_r$ is positive, so that from
  the fact that $c>0$ and $N^*>0$ we see that $s(A_r)$, the spectral bound of $A_r$, is a continuous
  monotone increasing function of $r$. By the assertion (3) of Lemma 3.3 we have
$$
  s(A_0)\leq -a<0 \quad \mbox{and} \quad \lim_{r\to+\infty}s(A_r)=+\infty.
$$
  Hence there exists a threshold value $r_*>0$ for $r$, such that $s(A_r)<0$ for $r<r_*$,
  $s(A_{r_*})=0$, and $s(A_r)\geq 0$ for $r>r_*$. Notice that when $a=r=0$, $A_r$ is exactly the
  operator $H$ given in the proof of Theorem 2.5. Since $s(H)=0$, by using Corollary 1.11 in
  Chapter VI of \cite{EN} we see that
$$
  \frac{a}{c\max_{x\in\overline{\Omega}}N^*(x)}\leq r_*\leq
  \frac{a}{c\min_{x\in\overline{\Omega}}N^*(x)}.
$$
  Notice that this implies that there exists $x_0\in\overline{\Omega}$ such that $a-r_*cN^*(x)=
  r_*c[N^*(x_0)-N^*(x)]$.
\medskip

  {\bf Theorem 3.5}\ \ {\em Let $r_*>0$ be as above. For $0<r<r_*$, the system of equations
  $(3.19)$ does not have a nontrivial nonnegative solution. Moreover, the solution $(I,R)$ of the
  system $(3.18)$ with nonnegative continuous initial data satisfies the following property:}
\begin{equation}
  \lim_{t\to\infty}I(x,t)=\lim_{t\to\infty}R(x,t)=0 \quad \mbox{uniformly for}\;\;
  x\in\overline{\Omega}.
\end{equation}

  {\em Proof}:\ \ The condition $0<r<r_*$ implies that $s(A_r)<0$, so that by the assertion (4) of
  Lemma 3.3 we have $R(0,A_r)\geq 0$. Notice that the equation (3.20) can be rewritten as follows:
$$
  I=-rR(0,A_r)[I(1+aL_b)I].
$$
  If $I$ is a nonnegative solution of this equation then the left-hand side is nonnegative while
  the right-hand side is non-positive, which can occur only in the case that both sides are
  vanishing. Hence if $0<r<r_*$ then the equation (3.20) does not have a nontrivial nonnegative
  solution and, consequently, the system of equation (3.19) does not have a nontrivial nonnegative
  solution. This proves the first assertion in Theorem 3.5.

  Next, we notice that from the first equation in (3.18) we see that for a nonnegative solution we
  have
$$
  \frac{\partial I(x,t)}{\partial t}\leq\int_{\Omega}v(x,y)I(y,t)\rmd y
  -[v_e(x)+a-rN(x,t)]I(x,t), \quad x\in\Omega,\;\; t>0,
$$
  or $I'(t)\leq A_rI(t)$ in abstract version. Since the $C_0$-semigroup generated by $A_r$ is
  positive, this implies that $I(t)\leq e^{tA_r}I_0$ for all $t\geq 0$, where $I_0$ is the initial
  data of $I$. It follows that
$$
  \max_{x\in\overline{\Omega}}I(x,t)\leq Ce^{ts(A_r)}\max_{x\in\overline{\Omega}}I_0(x),
   \quad t\geq 0,
$$
  where $C$ is a positive constant independent of $I_0$ and $t$. Since $s(A_r)<0$ for $0<r<r_*$,
  this implies that
$$
  \lim_{t\to\infty}I(x,t)=0 \quad \mbox{uniformly for}\;\; x\in\overline{\Omega}.
$$
  Using this relation and the second equation in (3.18), we can easily prove that also
$$
  \lim_{t\to\infty}R(x,t)=0 \quad \mbox{uniformly for}\;\; x\in\overline{\Omega}.
$$
  This completes the proof of Theorem 3.5. $\quad\Box$
\medskip

  {\bf Corollary 3.6}\ \ {\em Let $r_*>0$ be as above. For $0<r<r_*$, the solution $(S,I,R)$ of
  the system $(3.15)$ with nonnegative continuous initial data satisfies the following property:
\begin{equation*}
  \lim_{t\to\infty}S(x,t)=M_0N^*(x) \quad \mbox{uniformly for}\;\;
  x\in\overline{\Omega}, \quad \mbox{and}
\end{equation*}
\begin{equation*}
  \lim_{t\to\infty}I(x,t)=\lim_{t\to\infty}R(x,t)=0 \quad \mbox{uniformly for}\;\;
  x\in\overline{\Omega}.
\end{equation*}
  where $M_0=\displaystyle\int_{\Omega}[S_0(x)+I_0(x)+R_0(x)]\rmd x$. $\quad\Box$}
\medskip

  We conjecture that for $r>r_*$, the system (3.15) has a unique nontrivial nonnegative
  steady-state solution whose all three components are positive functions, and all nontrivial
  nonnegative time-dependent solutions of (3.15) converges to this steady-state solution as
  time goes to infinity. However, due to technical reasons we are unable to determine either
  this conjecture is true or false. A main difficult is caused by the nonlocal operator $L_b$ that
  appears in (3.20). In what follows we prove a weaker result.
\medskip

  {\bf Theorem 3.7}\ \ {\em There exists $r^*\geq r_*$ such that for $r>r^*$ the equation $(3.19)$
  has a nontrivial nonnegative solution $(I_r,R_r)$ with $I_r>0$, $R_r>0$ and
\begin{equation}
  \lim_{r\to\infty}I_r(x)=\frac{bM_0}{a+b}N^*(x) \quad \mbox{and} \quad
  \lim_{r\to\infty}R_r(x)=\frac{aM_0}{a+b}N^*(x) \quad
   \mbox{uniformly for}\;\; x\in\Omega,
\end{equation}
  where $M_0=\displaystyle\int_{\Omega}N_0(x)\rmd x=\int_{\Omega}[S_0(x)+I_0(x)+R_0(x)]\rmd x$.
  Moreover, for a possibly larger $r^*$ and any $r>r^*$, the solution $(I,R)$ of the system
  $(3.18)$ with initial value being small perturbation of $(I_r,R_r)$ in $C(\overline{\Omega})
  \times C(\overline{\Omega})$ satisfies the following property:
\begin{equation}
  \lim_{t\to\infty}I(x,t)=I_r(x) \quad \mbox{and} \quad
  \lim_{t\to\infty}R(x,t)=R_r(x) \quad \mbox{uniformly for}\;\;
  x\in\overline{\Omega}.
\end{equation}
}
\medskip

  {\em Proof}:\ \ Let $X=C(\overline{\Omega})$. Consider the mapping $F:X\times\mathbb{R}\to X$
  defined by
\begin{equation*}
  F(I,t)(x)=t\Big\{\int_{\Omega}v(x,y)I(y)\rmd y-[v_e(x)+a]I(x)\Big\}
  +[M_0N^*(x)-(aL_b+1)I(x)]I(x),  \quad x\in\Omega,
\end{equation*}
  where $I\in X$ and $t\in\mathbb{R}$. Let $I_{\infty}(x)=bM_0N^*(x)/(a+b)$. A simple
  computation shows that
$$
  F(I_{\infty},0)=0 \quad \mbox{and} \quad D_IF(I_{\infty},0)\varphi=-I_{\infty}(aL_b+1)\varphi
  \quad \mbox{and} \quad \varphi\in X.
$$
  The operator $D_IF(I_{\infty},0):X\to X$ has a bounded inverse. Indeed, since $L_b$ is the inverse
  of the operator
$$
  H_b:\varphi\mapsto [v_e(x)+b]\varphi(x)-\int_{\Omega}v(x,y)\varphi(y)\rmd y, \quad \forall\varphi\in X,
$$
  we see that for any $f\in X$, the equation $D_IF(I_{\infty},0)\varphi=f$ is equivalent to the
  following equation:
$$
  [v_e(x)+a+b]\varphi(x)-\int_{\Omega}v(x,y)\varphi(y)\rmd y=-H_b\Big(\frac{f}{I_{\infty}}\Big).
$$
  The unique solution of this equation is given by $\varphi=-L_{a+b}H_b(f/I_{\infty})$, and the
  operator $f\mapsto-L_{a+b}H_b(f/I_{\infty})$ in $X$ is clearly a bounded linear operator. Hence,
  by the inverse function theorem, there exist positive numbers $\varepsilon$ and $\delta$ such
  that for any $|t|<\varepsilon$, the equation $F(I,t)=0$ has a unique solution $I=I(t)$ in the
  $\delta$-neighborhood of $I_{\infty}$, i.e., $\|I(t)-I_{\infty}\|_X<\delta$, and $I(0)=
  I_{\infty}$. Letting $r^*=1/\varepsilon$ and $I_r=I(1/r)$, we see that for any $r>r^*$, the
  equation (3.20) has a solution $I_r\in X$ satisfying the property
$$
  \lim_{r\to\infty}I_r(x)=I_{\infty}(x)=\frac{bM_0}{a+b}N^*(x) \quad
   \mbox{uniformly for}\;\; x\in\Omega.
$$
  Since $N^*(x)>0$ for all $x\in\overline{\Omega}$, from the above relation we see that by choosing
  $r^*$ larger when necessary, we have $I_r(x)>0$ for all $x\in\overline{\Omega}$ and $r>r^*$. Now
  let
$$
  R_r=aL_bI_r.
$$
  Then clearly $R_r(x)>0$ for all $x\in\overline{\Omega}$ and, since $L_bN^*=(1/b)N^*$, we have
$$
  \lim_{r\to\infty}R_r(x)=\frac{aM_0}{a+b}N^*(x) \quad \mbox{uniformly for}\;\; x\in\Omega.
$$
  Moreover, from the relation between (3.20) and (3.19) we see that $(I,R)=(I_r,R_r)$ is a solution
  of (3.19). This proves the first part of Theorem 3.7.

  To prove the second part of Theorem 3.7, we let $\varphi(x,t)=I(x,t)-I_r(x)$ and $\psi(x,t)=
  R(x,t)-R_r(x)$, so that
$$
  I(x,t)=I_r(x)+\varphi(x,t), \quad  R(x,t)=R_r(x)+\psi(x,t).
$$
  Substituting these expressions and (2.9) into (3.18) and using the fact that $(I_r,R_r)$ is a
  solution of (3.19), we see that $(\varphi,\psi)$ satisfies the following system of equations:
\begin{equation*}
\left\{
\begin{array}{l}
  \displaystyle\frac{\partial\varphi(x,t)}{\partial t}=\int_{\Omega}v(x,y)\varphi(y,t)\rmd y
  -[v_e(x)+a-rM_0N^*(x)+2rI_r(x)+rR_r(x)]\varphi(x,t)
\\
   \quad\quad\quad
  -rI_r(x)\psi(x,t)-r\varphi^2(x,t)-r\varphi(x,t)\psi(x,t)+r[I_r(x)+\varphi(x,t)]h(x,t),
  \quad x\in\Omega,\;\; t>0,
\\ [0.3cm]
  \displaystyle\frac{\partial\psi(x,t)}{\partial t}=\int_{\Omega}v(x,y)\psi(y,t)\rmd y
  -[v_e(x)+b]\psi(x,t)+a\varphi(x,t), \quad x\in\Omega,\;\; t>0.
\end{array}
\right.
\end{equation*}
  To prove (3.31) we only need to prove that there exists $\mu>0$ such that the solution of the
  above system satisfies the following estimates:
\begin{equation}
  \|\varphi(\cdot,t)\|_{\infty}+\|\psi(\cdot,t)\|_{\infty}\leq
  Ce^{-\mu t}[\|\varphi(\cdot,0)\|_{\infty}+\|\psi(\cdot,0)\|_{\infty}] \quad
  \mbox{for}\;\; t\geq 0
\end{equation}
  if $\|\varphi(\cdot,0)\|_{\infty}+\|\psi(\cdot,0)\|_{\infty}$ is sufficient small.
  Let $\bfA:X\times X\to X\times X$ be the following bounded linear operator: For $(\varphi,\psi)
  \in X\times X$, $\bfA(\varphi,\psi)=(A_1(\varphi,\psi),A_2(\varphi,\psi))$, where
$$
\begin{array}{l}
  A_1(\varphi,\psi)=\displaystyle\int_{\Omega}v(x,y)\varphi(y)\rmd y-[v_e(x)+a
  -rM_0N^*(x)+2rI_r(x)+rR_r(x)]\varphi(x)-rI_r(x)\psi(x),
\\ [0.3cm]
  A_2(\varphi,\psi)=\displaystyle\int_{\Omega}v(x,y)\psi(y)\rmd y-[v_e(x)+b]\psi(x,t)+a\varphi(x).
\end{array}
$$
  Since we already know that $\|h(\cdot,t)\|_{\infty}\leq Ce^{-\kappa t}\|h(\cdot,0)\|_{\infty}$
  for $t\geq 0$, to prove (3.31) we only need to prove that $s(\bfA)<0$. In what follows we prove
  that for a possibly larger $r^*$, if $r>r^*$ then for all $\lambda\in\mathbb{C}$ with
  $\mbox{\rm Re}\lambda>-\min\{a,b\}$ we have $\lambda\in\rho(\bfA)$, so that $s(\bfA)\leq
  -\min\{a,b\}$.

  Given $(f,g)\in X\times X$, the equation $(\lambda I-\bfA)(\varphi,\psi)=(f,g)$ can be explicitly
  rewritten as follows:
\begin{equation}
\left\{
\begin{array}{l}
  \displaystyle[\lambda+a+v_e(x)+2rI_r(x)+rR_r(x)-rM_0N^*(x)]\varphi(x)
  -\int_{\Omega}v(x,y)\varphi(y)\rmd y+rI_r(x)\psi(x)=f(x),
\\ [0.3cm]
  \displaystyle[\lambda+b+v_e(x)]\psi(x,t)-\int_{\Omega}v(x,y)\psi(y)\rmd y-a\varphi(x)=g(x).
\end{array}
\right.
\end{equation}
  From (3.31) we see that if $r^*$ is sufficiently larger and $r>r^*$ then $2I_r(x)+R_r(x)
  -M_0N^*(x)\geq C_0$ for some constant $C_0>0$. Later on we assume that this is the case. From the
  second equation in (3.31) we see that if $\mbox{\rm Re}\lambda>-b$ then
\begin{equation}
  \psi=aL_{\lambda+b}\varphi+L_{\lambda+b}g.
\end{equation}
  Notice that there exists positive constant $C$ independent of $\lambda$ such that for all
  $\mbox{\rm Re}\lambda>-b$,
\begin{equation}
  \|L_{\lambda+b}\|_{\mathscr{L}(X)}\leq\frac{C}{\mbox{\rm Re}\lambda+b}.
\end{equation}
  Substituting the expression (3.31) into the first equation in (3.31) we see that unique
  solvability of the system (3.31) is equivalent to unique solvability of the following equation:
\begin{equation}
\begin{array}{l}
  [\lambda+a+v_e(x)+2rI_r(x)+rR_r(x)-rM_0N^*(x)]\varphi(x)
\\
  \displaystyle -\int_{\Omega}v(x,y)\varphi(y)\rmd y+arI_r(x)L_{\lambda+b}\varphi(x)
  =f(x)-rL_{\lambda+b}g(x).
\end{array}
\end{equation}
  For any $\mbox{\rm Re}\lambda>-a$ and $r>r^*$, since $2I_r(x)+R_r(x)-M_0N^*(x)\geq C_0>0$, a
  similar consideration as before shows that the operator $\varphi\mapsto\displaystyle [\lambda+a+
  v_e(x)+2rI_r(x)+rR_r(x)-rM_0N^*(x)]\varphi(x)-\int_{\Omega}v(x,y)\varphi(y)\rmd y$ in $X$ is
  invertible, and the norm of the inverse is bounded by $C(\mbox{\rm Re}\lambda+a+C_0r)^{-1}$,
  where $C$ is a positive constant independent of $\lambda$ and $r$. Let $B_{\lambda,r}$ be this
  inverse. Then (3.31) is equivalent to the following equation:
\begin{equation}
  \varphi(x)+arB_{\lambda,r}[I_rL_{\lambda+b}\varphi](x)
  =B_{\lambda,r}f(x)-rB_{\lambda,r}L_{\lambda+b}g(x).
\end{equation}
  Using (3.31) and the fact that $\|B_{\lambda,r}\|_{\mathscr{L}(X)}\leq C(\mbox{\rm Re}\lambda+a
  +C_0r)^{-1}$ we easily see that there exists $\Lambda>\max\{a,b\}$ such that if $|\lambda|>
  \Lambda$ then for any $r>r^*$,
$$
  \|arB_{\lambda,r}[I_rL_{\lambda+b}\varphi]\|_X\leq\frac{1}{2}\|\varphi]\|_X.
$$
  Hence, by using the contraction mapping theorem we see that if $|\lambda|>\Lambda$ then for any
  $r>r^*$, the equation (3.31) is uniquely solvable and, consequently, the system (3.31) is
  uniquely solvable. Next we consider the case $|\lambda|\leq\Lambda$, $\mbox{\rm Re}\lambda>
  -\min\{a,b\}$. Since the operator
$$
  \varphi\mapsto [\lambda+a+v_e(x)]\varphi(x)-\int_{\Omega}v(x,y)\varphi(y)\rmd y
$$
  is uniformly bounded in $X$ for $\lambda$ in this region, it can be easily seen that if $r^*$ is
  sufficiently larger and $r>r^*$ then from the first equation in (3.31) we have
\begin{equation}
  \varphi(x)=-m_r(x)\psi(x)+\tilde{B}_{\lambda,r}\psi(x)+\hat{B}_{\lambda,r}f(x),
\end{equation}
  where $m_r(x)=[2I_r(x)+R_r(x)-M_0N^*(x)]^{-1}I_r(x)$, $\tilde{B}_{\lambda,r}$ and
  $\hat{B}_{\lambda,r}$ are bounded linear operators in $X$, and
\begin{equation}
  \|\tilde{B}_{\lambda,r}\|_{\mathscr{L}(X)}\leq Cr^{-1},
\end{equation}
  where $C$ is a constant independent of $\lambda$ and $r$. Notice that $m_r>0$, which implies
  that for any $\lambda\in\mathbb{C}$ such that $\mbox{\rm Re}\lambda>-b$, the operator
$$
  \psi\mapsto [\lambda+b+v_e(x)+m_r(x)]\psi(x)-\int_{\Omega}v(x,y)\psi(y)\rmd y
$$
  has bounded inverse in $X$. Substitute (2.32) into the second equation in (3.31) and using (3.31)
  and a contraction mapping argument, we see that by choosing $r^*$ larger when necessary, for any
  $r>r^*$ the equation for $\psi$ obtained in this way is uniquely solvable. Hence for $\lambda\in
  \mathbb{C}$ satisfying $\mbox{\rm Re}\lambda>-\min\{a,b\}$, $|\lambda|\leq\Lambda$ and for $r>
  r^*$, the system (3.31) is also uniquely solvable. This proves the desired assertion and
  completes the proof of Theorem 3.7. $\quad\Box$
\medskip

  {\bf Corollary 3.8}\ \ {\em Let $r^*>0$ be as above. For $r>r_*$, the system $(3.15)$ has a
  nontrivial nonnegative steady-state solution $(S_r,I_r,R_r)$ such that $S_r>0$, $I_r>0$ and
  $R_r>0$. Moreover, any time-dependent solution $(S,I,R)$ of the system $(3.15)$ starting from a
  small neighborhood of $(S_r,I_r,R_r)$ in supremum norm sense satisfies the following property:
\begin{equation*}
  \lim_{t\to\infty}S(x,t)=S_r(x), \quad \lim_{t\to\infty}I(x,t)=I_r(x), \quad
   \lim_{t\to\infty}R(x,t)=R_r(x)
\end{equation*}
  exponentially uniformly for $x\in\overline{\Omega}$.}
\medskip

  {\em Proof}:\ \ Let $(I_r,R_r)$ be the nontrivial nonnegative solution of the system (3.19)
  ensured by Theorem 3.7, and consider the equation
\begin{equation*}
  \int_{\Omega}v(x,y)S_r(y)\rmd y-v_e(x)S_r(x)-rS_r(x)I_r(x)+bR_r(x)=0, \quad x\in\Omega,
\end{equation*}
  or
\begin{equation*}
  [v_e(x)+rI_r(x)]S_r(x)-\int_{\Omega}v(x,y)S_r(y)\rmd y=bR_r(x), \quad x\in\Omega,
\end{equation*}
  Since $I_r>0$ and $R_r>0$, it is easy to see that this equation has a unique solution which is
  strictly positive in $\overline{\Omega}$. This proves the first assertion of Corollary 3.8. The
  proof of the second assertion is easy and we omit it here. $\quad\Box$
\medskip

\section{Conclusions}
\setcounter{equation}{0}

\hskip 2em
  We have studied two classes of mathematical models: population migration models and migration
  epidemics models. By population migration model we mean the mathematical model describing
  evolution of a population system in which spatial movement of individuals depends only on the
  departure and arrival locations and does not have apparent dependence on population density (in
  contrast to diffusion population models in which spatial movement of individuals is determined
  by gradient of the population density). Mathematical formulation of a such model is a linear
  differential-integral equation. We mainly consider the ergodic population system as defined in
  Definition 2.1. In the proliferation-stationary case, by using the theories of irreducible
  positive operators and irreducible positive semigroups, we prove that in addition to the obvious
  assertion that total population keeps constant, the population will finally distribute to the
  whole habitat domain in a fixed proportion $N^*$ uniquely determined by the migration rate function,
  regardless how the population is initially distributed in the habitat domain. More precisely,
  the population density $N$ has the following asymptotic expression:
\begin{equation}
  N(x,t)=M_0N^*(x)[1+o(1)] \quad \mbox{for}\;\; x\in\Omega, \;\, \mbox{as}\;\;  t\to\infty,
\end{equation}
  where $N^*$ is a strictly positive function uniquely determined by the migration rate function,
  and $M_0$ is the initial total population. In the case that the proliferation rate is not
  identically zero, the result depends on a parameter $s(H)$, the spectral bound of the operator
  $H$ (see (2.30) for the expression of the operator $H$): If $s(H)<0$ then the population will
  finally become extinct, if $s(H)=0$ then under the additional condition $s(H)>\displaystyle
  \max_{x\in\overline{\Omega}}[r(x)-v_e(x)]$, where $r$ is the proliferation rate function and
  $v_e$ is the emigration rate function, a similar assertion as in the proliferation-stationary
  case holds, and if $s(H)>0$ then under the same condition as given above the so-called balanced
  exponential growth phenomenon occurs, i.e., the population density $N$ has the following
  asymptotic expression (see (2.35)):
\begin{equation*}
  N(x,t)=M_0N^*(x)e^{s(H)\, t}[1+o(1)] \quad \mbox{for}\;\; x\in\Omega, \;\,
  \mbox{as}\;\;  t\to\infty,
\end{equation*}
  where $N^*$ is a similar function as above and $M_0$ is a constant related to the initial
  population density $N_0$. The following estimates for $s(H)$ are also proved:
\begin{equation*}
  \max_{x\in\overline{\Omega}}[r(x)-v_e(x)]\leq s(H)\leq\max_{x\in\overline{\Omega}}r(x).
\end{equation*}
  Besides, we have also discussed two examples of non-ergodic situations. These examples show
  that non-ergodic population systems have different dynamical behavior.

  We note that in some existing literatures some nonlinear version of the population migration
  model studied here have been studied by using upper and lower solutions method, cf. \cite{Bat,
  Can2,Can3,Cos4,Hut2,Cov1,Cov2,Cov3,Cov4} and references cited therein. The results obtained here
  are clearly different from those obtained in the references.

  Based on the above analysis on population migration models, we have also studied three simple
  mathematical models describing spread of epidemics in an ergodic migration population system.
  These models include: The migration SI model (3.1), the migration SIR model (3.9), and the
  migration SIRE model (3.15). For the migration SI model (3.1) we have proved the following result:
\begin{equation}
  \lim_{t\to\infty}(S(x,t),I(x,t))=(0,M_0N^*(x)) \quad \mbox{uniformly for}\;\; x\in\Omega,
\end{equation}
  where $S$ and $I$ denote the densities of the uninfected (or susceptible) and infected populations,
  respectively, and $N^*$, $M_0$ are as in (4.1). This means that in an ergodic population system,
  when an outbreak of an epidemics takes place all individuals will finally become infected,
  regardless where they live and initially how and where the infected individuals are distributed.
  Thus, in an ergodic population system as is the modern global human population system, when an
  outbreak of epidemics occurs, to avoid uncontrollable spread of the epidemics, segregation of
  infected individuals from uninfected ones is a necessary measure. For the migration SIR model
  (3.9) we have proved the following result:
\begin{equation*}
  \lim_{t\to\infty}(S(\cdot,t),I(\cdot,t),R(\cdot,t))=(c_1N^*,0,c_2N^*)
  \quad \mbox{in $L^1$ sense},
\end{equation*}
  where $S$, $I$ and $R$ denote the densities of the uninfected (or susceptible), the infectious,
  and the recovered and immune populations, respectively, $N^*$ is as in (4.1), and $c_1$, $c_2$
  are positive constants, $c_1+c_2=M_0$. This means that if a infectious disease is not so harmful
  in the sense that it is not fatal (i.e., infection of the disease does not cause death of
  population) and infected individuals can recover and become immune of it, then when it occurs in
  an ergodic population system, all infected individuals will finally recover, and there is a
  portion of the population keeping uninfected all the time. Moreover, ergodicity ensures that the
  uninfected and recovered populations will both distribute to the whole habitat domain in some
  proportions $\lambda_1N^*$ and $\lambda_2N^*$, respectively, where $\lambda_1$, $\lambda_2$ are
  positive constants depending only on the initial total population $M_0$ while independent of the
  initial distribution of any of the three groups of population, such that $\lambda_1+\lambda_2=1$.

  Finally, the migration SIRE model (3.15) describes the evolution of also three groups of
  individuals in a population system inhabiting a spatial domain $\Omega$: the susceptible group,
  the infectious group, and the recovered group. Unlike the migration SIR model mentioned above,
  in the SIRE model individuals in the recovered group are not completely immune to the disease;
  they can be infected with the disease again in a suitable percentage or after a certain period
  of recovery, so that the susceptible group includes both individuals who are never infected and
  those who have recovered from infection but have lost immunity. Due to mathematical difficulties
  our results of analysis to this model is not very satisfactory: We have only proved that there
  exist constants $0<r_*\leq r^*<\infty$ such that if the infection rate $0<r<r_*$ then
\begin{equation*}
  \lim_{t\to\infty}(S(x,t),I(x,t),R(x,t))=(M_0N^*(x),0,0)
  \quad \mbox{uniformly for}\;\; x\in\Omega,
\end{equation*}
  where $S$ and $I$ denote the densities of the susceptible, the infected and the recovered
  populations, respectively, and $N^*$, $M_0$ are as in (4.1), i.e., if the infection rate is
  sufficiently small then the population system will finally become free of the epidemics, and if
  $r>r^*$ then there exists a stationary state in which all the three groups of population has
  non-identically vanishing densities $(S_r(x),I_r(x),R_r(x))$, $S_r(x)+I_r(x)+R_r(x)=M_0N^*(x)$,
  such that
\begin{equation*}
  \lim_{t\to\infty}(S(x,t),I(x,t),R(x,t))=(S_r(x),I_r(x),R_r(x))
  \quad \mbox{uniformly for}\;\; x\in\Omega,
\end{equation*}
  if the initial population densities $S_0(x),I_0(x),R_0(x)$ are sufficiently closed to $S_r(x),
  I_r(x),R_r(x)$, respectively.

  Despite that the analysis made here is very elementary, we have still acquired much knowledge on
  dynamical behavior of a population system in which individuals can migrate from one place to
  another, and also obtained some useful information on the effects of population migration to
  spread of epidemics. Limited to spaces, we have to stop the discussion here and leave more
  penetrating investigation for future work.
\medskip

  {\bf Acknowledgement}\ \ It is our great pleasure to acknowledge sincere thanks to Professor
  Chris Cosner for his kindness of providing with the references \cite{Cos3, Cos4}, which greatly
  helped us in tracking previous literatures on the topic studied in this manuscript.

\clearpage
\begin{table}
\begin{center}
\begin{tabular}{|c|c|c|c|c|c|c|c|}
\hline
{\bf\small Province} & {\small Beijing} & {\small Tianjin} & {\small Hebei} & {\small Shanxi} &
  {\small ${\mbox{Inner}\atop\mbox{Mongolia}}$} & {\small Liaoning} & {\small Jilin} \\
\hline
{\bf\small Infection} & 2434 & 176 & 210 & 445 & 289 & 3 & 35 \\
\hline
{\bf\small Death} & 147 & 12 & 10 & 25 & 25 & 0 & 3 \\
\hline
\hline
{\bf\small Province} & {\small Shanghai} & {\small Jiangsu} & {\small Zhejiang} & {\small Anhui} &
  {\small Fujian} & {\small Jiangxi} & {\small Shandong} \\
\hline
{\bf\small Infection} & 7 & 7 & 4 & 10 & 3 & 1 & 1 \\
\hline
{\bf\small Death} & 2 & 1 & 0 & 0 & 0 & 0 & 0 \\
\hline
\hline
{\bf\small Province} & {\small Guangdong} & {\small Guangxi} & {\small Hunan} & {\small Hubei} &
  {\small Henan} & {\small Chongqing} & {\small Sichuan} \\
\hline
{\bf\small Infection} & 1514 & 22 & 6 & 6 & 15 & 3 & 17 \\
\hline
{\bf\small Death} & 56 & 3 & 1 & 0 & 0 & 0 & 2 \\
\hline
\hline
{\bf\small Province} & {\small Shannxi} & {\small Gansu} & {\small Ningxia} & {\small Taiwan} &
  {\small Hongkong} & {\small Macao} & {\small ${\mbox{Other}\atop\mbox{Provinces}}$} \\
\hline
{\bf\small Infection} & 12 & 8 & 6 & 307 & 1755 & 1 & 0 \\
\hline
{\bf\small Death} & 0 & 1 & 1 & 47 & 299 & 0 & 0 \\
\hline
\end{tabular}
\end{center}
\caption{Statistics of SARS cases in China}
\end{table}
\begin{table}
\begin{center}
\begin{tabular}{|c|c|c|c|c|c|c|}
\hline
{\bf\small Country} & {\small Australia} & {\small Brazil} & {\small Canada} & {\small China}
  & {\small Columbia} & {\small Finland} \\
\hline
{\bf\small Infection} & 5 & 1 & 250 & 6761 & 1 & 1  \\
\hline
{\bf\small Death} & 0 & 0 & 38 & 570 & 0 & 0  \\
\hline
\hline
{\bf\small Country} & {\small France} & {\small Germany} & {\small India} & {\small Indonesia}
  & {\small Italy} & {\small Kuwait} \\
\hline
{\bf\small Infection} & 7 & 10 & 3 & 2 & 4 & 1  \\
\hline
{\bf\small Death} & 1 & 0 & 0 & 0 & 0  & 0\\
\hline
\hline
{\bf\small Country} & {\small Malaysia} & {\small Mongolia} & {\small New Zealand}
  & {\small Philippine} & {\small ${\mbox{Republic of}\atop\mbox{Ireland}}$}
  & {\small ${\mbox{Republic of}\atop\mbox{Korea}}$} \\
\hline
{\bf\small Infection} & 5 & 9 & 1 & 14 & 1 & 3 \\
\hline
{\bf\small Death} & 2 & 0 & 0 & 2 & 0 & 0 \\
\hline
\hline
{\bf\small Country} & {\small Romania} & {\small ${\mbox{Russian}\atop\mbox{Federation}}$}
  & {\small Singapore} & {\small South Africa} & {\small Spain} & {\small Sweden} \\
\hline
{\bf\small Infection} & 1 & 1 & 206 & 1 & 1 & 3 \\
\hline
{\bf\small Death} & 0 & 0 & 32 & 1 & 0 & 0  \\
\hline
\hline
{\bf\small Country} & {\small Switzerland} & {\small Thailand} & {\small ${\mbox{United}\atop\mbox{Kingdom}}$}
  & {\small United States} & {\small Viet Nam} & {\small ${\mbox{Other}\atop\mbox{Countries}}$} \\
\hline
{\bf\small Infection} & 1 & 9 & 4 & 71 & 63 & 0 \\
\hline
{\bf\small Death} & 0 & 2 & 0 & 0 & 5 & 0  \\
\hline
\end{tabular}
\end{center}
\caption{Statistics of SARS cases in World}
\end{table}

\end{document}